\documentclass[12pt]{article}
\usepackage{amsmath, amsthm, amssymb, amscd, graphicx}
\pdfoutput=1 

\usepackage{geometry} % to change the page dimensions

\title{Discontinuity growth of interval exchange maps}
\author{Christopher F. Novak}
\newtheorem{thm}{Theorem}[section]
\newtheorem{cor}[thm]{Corollary}
\newtheorem{lem}[thm]{Lemma}
\newtheorem{prop}[thm]{Proposition}
\newtheorem*{thm_nonum}{Theorem}

\theoremstyle{remark}
\newtheorem{rem}[thm]{Remark}

\theoremstyle{definition}
\newtheorem{defn}[thm]{Definition}

\begin{document}

\maketitle
%\tableofcontents

\begin{abstract}
For an interval exchange map $f$, the number of discontinuities $d(f^n)$ 
either exhibits linear growth or is bounded independently of $n.$  
This dichotomy is used to prove that the group $\mathcal{E}$ of interval exchanges does not contain distortion elements, giving
examples of groups that do not act faithfully via interval exchanges. As a further application of 
this dichotomy, a classification of centralizers in $\mathcal{E}$ is given. This classification is used to show
 that $\text{Aut}(\mathcal{E}) \cong \mathcal{E} \rtimes \mathbb{Z}/ 2 \mathbb{Z}.$ 
 \end{abstract}

%  and a complete characterization of maps with bounded discontinuity growth is given.
%  
%   include this 

\section{Introduction}
\begin{figure}[htb] \label{interval exchange}
\begin{center}
\includegraphics{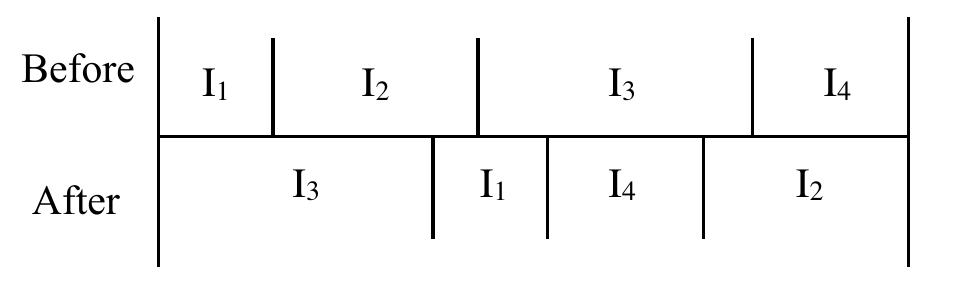}
\caption{An interval exchange with $\pi = (2, 4, 1, 3)$}
\end{center}
\end{figure}

An interval exchange transformation is a map $\mathbb{T}^1 \rightarrow \mathbb{T}^1$   
defined by a partition of the circle into a finite union of half-open intervals and 
 a rearrangement of these intervals by translation.        
See Figure \ref{interval exchange} for a graphical example.

The dynamics of interval exchanges were first studied in the late seventies by
 Keane{\bf \cite{Keane 75}},{\bf\cite{Keane 77}},  Rauzy{\bf \cite{Rauzy 79}}, Veech{\bf \cite{Veech 78}}, and others. 
 This initial stage of research culminated in the independent proofs by Masur{\bf \cite{Masur 82}} and Veech{\bf\cite{Veech 82}} that
 almost every interval exchange is uniquely ergodic.
 See the recent survey of Viana{\bf \cite{Viana 06}} for a unified presentation of these results. 
  The current study of interval exchanges is closely related to dynamics 
 on the moduli space of translation surfaces; an introduction to this topic and its connection to interval exchanges is found in a survey of
 Zorich{\bf\cite{Zorich 06}}.

To precisely define an interval exchange, let $\pi \in \Sigma_n$ be a permutation and let 
$\lambda = (\lambda_1, \ldots, \lambda_n)$ be a vector in the simplex 
\[ \Lambda_n = \lbrace(\lambda_1, \ldots, \lambda_n) \in \mathbb{R}^n: \lambda_i >0 \text{ and } \sum \lambda_i = 1 \rbrace.   \]
The vector $\lambda$ induces a partition of $\mathbb{T}^1 \cong [0,1)$ into intervals of length $\lambda_j$: 
\[     I_j = [\beta_{j-1}, \beta_j) \quad  1\leq j \leq n,  \]
\[  \text{where }   \beta_0 = 0, \quad \beta_j = \sum_{i=1}^j \lambda_i \text{ for } 1\leq j \leq n.\]
The interval exchange $f_{(\pi, \lambda)}$ reorders the $I_j$ by translation, such that their indicies  are ordered 
by $\pi^{-1}(1), \pi^{-1}(2), \ldots, \pi^{-1}(n).$ 
Consequently,  $f_{(\pi, \lambda)}$ is defined by the formula 
\[  f_{(\pi, \lambda)} (x) =  x -\left( \sum_{i < j} \lambda_i \right) + \left( \sum_{\pi(i) < \pi(j)} \lambda_i \right)
 = x + \omega_j, \text{ for }x\in I_j. \] 
The vector $\omega(f) = (\omega_1, \ldots, \omega_n)$ is called the translation vector of $f_{(\pi, \lambda)}$. 

%By the above formula, 
%it is seen that there is an antisymmetric linear map $\Omega_{\pi}: \mathbb{R}^n \rightarrow \mathbb{R}^n$, depending 
%only on the permutation $\pi$,  such that $\Omega_{\pi} (\lambda) = \omega$.

 Let $d(f)$ denote the number of discontinuity points of $f\!:\!\mathbb{T}^1 \rightarrow \mathbb{T}^1$, where $\mathbb{T}^1$ is endowed 
 with its standard topology.  
  If $d(f) = k$, then it is easy to see that 
 for iterates $f^n,$ the discontinuity number $d(f^n)$ is bounded above by $k|n|$. It is possible for $d(f^n)$ to have linear growth at 
 a rate which is strictly less than the maximum $d(f)$. For example, the map in Figure 2
  has three discontinuities, but iteration 
 will suggest that $d(f^n) \sim 2n$.  
  Additionally, it is possible for $d(f^n)$ to be bounded independently of $n$. For example, a restricted rotation
 $r_{\alpha, \beta}$, as defined by Figure 3, satisfies $d(r_{\alpha, \beta}^n) \leq 3$ for all $n \in \mathbb{Z}$. The key result of this
 paper is the observation that no intermediate growth rate may occur.
 
 \begin{figure}[htb]
\begin{center}
\includegraphics{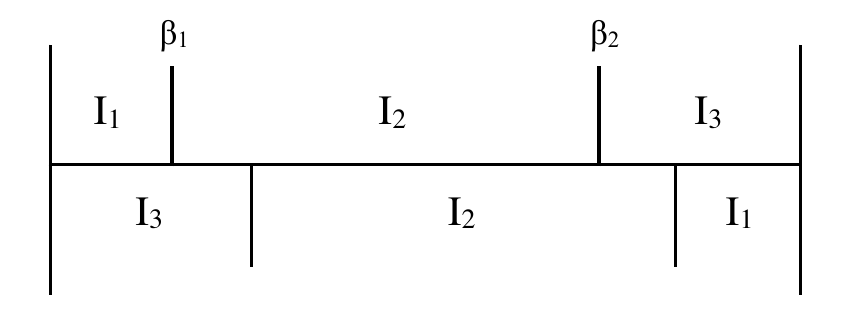}
\caption{$d(f^n)$ exhibits linear growth  ($\beta_1$ and $\beta_2$ independent over $\mathbb{Q}$)}
\end{center}
 \label{fig: linear growth example}
\end{figure}

 \begin{figure}[htb]
\begin{center}
\includegraphics{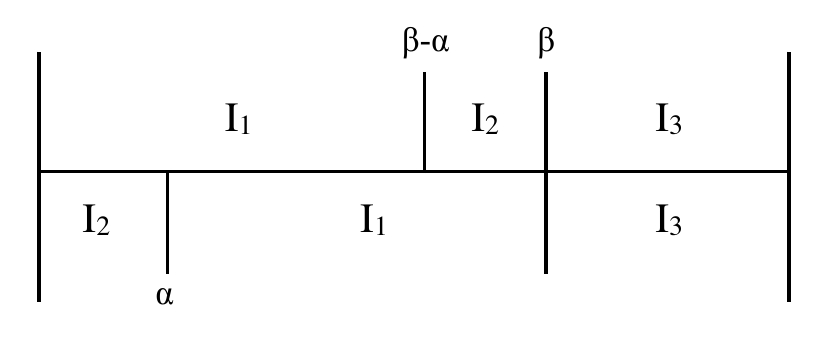}
\caption{The restricted rotation $r_{\alpha, \beta}:$  $\pi = (2, 1, 3)$,  $\lambda = (\beta - \alpha, \alpha, 1- \beta)$} 
\end{center}
 \label{fig: restricted rotation}
\end{figure}

 \begin{thm}
  \label{discontinuity growth dichotomy: basic version}
For any interval exchange $f,$ either $d(f^n)$ exhibits linear growth or $d(f^n)$ is bounded independently of $n$.
 \end{thm}

This theorem is a simplified statement of Proposition \ref{discontinuity growth dichotomy: detailed version}. The 
linear growth case is generic; for instance, 
 given any irreducible permutation $\pi$ which permutes three or more intervals, $f_{(\pi, \lambda)}$ has linear discontinuity growth
if the boundary points between intervals satisfy the \emph{infinite distinct orbit condition }{\bf \cite{Keane 75}}. This condition is satisfied for 
the full measure set of $\lambda \in \Lambda_n$ that have rationally independent 
partition lengths $\{\lambda_i\}.$ 

This raises the question of what may be said about the interval exchanges $f$ for which 
$d(f^n)$ is bounded. A result of Li {\bf \cite{Li 99}} stated in Section 3
 asserts that under certain additional conditions, the only such topologically minimal examples are maps 
 conjugate to an irrational rotation. 
 
By only assuming that $d(f^n)$ is bounded, it is still possible to give a complete description of the interval exchanges with bounded discontinuity growth.   
For $\gamma \in \mathbb{R}/\mathbb{Z} \cong [0, 1)$, let $r_\gamma$ denote the
 rotation $x \mapsto x + \gamma$, which is represented by the data $\pi = (2,1), \lambda = (1-\gamma, \gamma)$. 
An interval exchange is a \emph{restricted rotation} if it is 
conjugate by some $r_\gamma$ to some $r_{\alpha, \beta}.$ The \emph{support} of an interval exchange $f$ is
the complement of the set of fixed points $\text{Fix}(f).$ The following classification of interval exchanges with bounded discontinuity growth is 
proved in Section 3.  

 \begin{thm}  \label{bounded growth classification}
 Let $f$ be an infinite-order interval exchange transformation for which $d(f^n)$ is bounded. Then for some $k\geq 1,$ $f^k$ is 
 conjugate to a product of infinite-order restricted rotations with pairwise disjoint supports.
 \end{thm}
 
  \begin{figure}[htb]
\begin{center}
\includegraphics{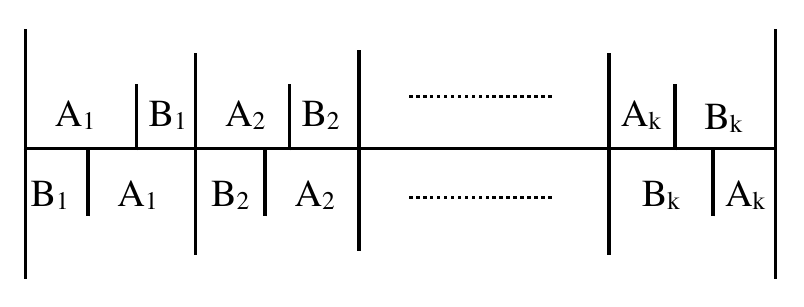}
\caption{A product of restricted rotations} 
\end{center}
 \label{fig: product of restricted rotation}
\end{figure}

The discontinuity growth dichotomy and the subsequent classification of maps with bounded discontinuity growth may be applied 
to the study of interval exchange group actions. The set $\mathcal{E}$ of all interval exchange transformations
on $\mathbb{T}^1$ forms a group under composition. For a group $G$, an \emph{interval exchange action} of $G$ is a 
group homomorphism $G \rightarrow \mathcal{E}.$ Such an action is \emph{faithful} if this homomorphism is injective, in which case
the image is a subgroup of $\mathcal{E}$ isomorphic to $G$.  

  Let $G$ be a finitely generated group, and let $S = \{g_1, \ldots, g_n\}$ be a set of generators.
   An element $f\in G$ is a \emph{distortion element} if $f$ has infinite order and
\[  \liminf_{n\rightarrow \infty}  \frac{|f^n|_S}{n} = 0, \]
where $|\cdot|_S$ denotes the minimal word length in terms of the generators and their inverses. 
For example, the central elements of the discrete Heisenberg group are distortion elements. In general, if $G$ is not finitely generated, an element
$f\in G$ is said to be a distortion element in $G$ if it is a distortion element in some finitely generated subgroup of $G$.

  \begin{thm}  \label{no distortion elements}
 The group $\mathcal{E}$ contains no distortion elements.
 \end{thm}

A proof of this result is given in Section 4.
The main consequence of this theorem is that any group $G$ containing a distortion element has no faithful interval exchange actions.  
 A particularly interesting case is the following, which is analogous to a result of Witte {\bf \cite{Witte 94}} for group actions
$SL(n, \mathbb{Z}) \rightarrow \text{Homeo}_{+}(S^1)$ for $n \geq 3$.    

\begin{cor}  \label{various lattices do not act}
Suppose $\Gamma$ is a non-uniform irreducible lattice in a semisimple Lie group $G$ with $\mathbb{R}$-rank $\geq 2$. Suppose further
that $G$ is connected, with finite center and no nontrivial compact factors. 
Then any interval exchange action $\Gamma \rightarrow \mathcal{E}$ has finite image. 
\end{cor} 

For example, the lattices $SL(n, \mathbb{Z}), n \geq 3,$ satisfy the above hypotheses; consequently, they do not act faithfully via interval exchange maps.
This corollary follows from a theorem of Lubotzky, Moses, and Raghunathan {\bf\cite{Lubostzky et al 00}} which states that lattices 
satisfying the above conditions contain distortion elements (in fact, elements with logarithmic word growth) and a theorem of Margulis 
{\bf\cite{Margulis 91}} which states that any irreducible lattice in a semisimple Lie group of $\mathbb{R}$-rank $\geq 2$ is almost simple; i.e., 
any normal subgroup of such a lattice is finite or has finite index.   \\

A further application of the discontinuity growth dichotomy is the complete classification of centralizers in the group $\mathcal{E}, $ which is
developed in Section 5. 
 This classification relies on analyzing the centralizer $C(f)$ in 
three cases that are distinguished by dynamical characteristics.

\begin{prop} Let $f$ be an interval exchange transformation.
\begin{itemize}
\item[ (i) ] $f$ has periodic points if and only if $C(f)$ contains a subgroup isomorphic to $\mathcal{E}$.
\item[ (ii) ] If $f$ is minimal and $d(f^n)$ is bounded, then $C(f)$ is virtually abelian and contains a 
	subgroup isomorphic to $\mathbb{R}/ \mathbb{Z}$.
\item[ (iii) ] If $f$ is minimal and $d(f^n)$ has linear growth, then $C(f)$ is virtually cyclic. 
\end{itemize}
\end{prop}

Minimality here refers to topological minimality: every orbit of $f$ is dense in $\mathbb{T}^1$. 
%The periodic case is characterized by the fact that the centralizer $C(f)$ contains a subgroup isomorphic to 
% $\mathcal{E}$. In the bounded growth case, $C(f)$ is virtually abelian and contains a subgroup
%isomorphic to $\mathbb{R}/\mathbb{Z}$. Finally, if $f$ is minimal and 
% has linear discontinuity growth, then $C(f)$ is virtually cyclic. 
The three parts of this result are restated and proved separately as Corollary \ref{periodic points iff big centralizer}, 
Proposition \ref{centralizer of minimal and bounded disc growth has rotation subgroup}, 
and Proposition \ref{minimality and linear discontinuity growth implies virtually cyclic centralizer}, respectively.
 These cases may be combined to give a general description 
 of centralizers in $\mathcal{E}$; this is stated and proved as Theorem \ref{general centralizer form}.

The classification of centralizers in $\mathcal{E}$ may be used to investigate the automorphism group $\text{Aut}(\mathcal{E})$. 
 Since $\mathcal{E}$ has trivial center, the inner automorphism group is isomorphic to $\mathcal{E}$.
 A further automorphism is induced by switching the orientation 
 of the circle $\mathbb{T}^1$. 
  More precisely, let $T\!:\! \mathbb{T}^1 \rightarrow \mathbb{T}^1$ be defined by $T(x) = -x$. 
 For $f\in \mathcal{E},\, T^{-1}fT$ is still an invertible piecewise translation, but it is now continuous from the left. 
 Let $\Psi_T$ be the automorphism of $\mathcal{E}$ defined by conjugation with $T$ followed by the natural isomorphism
 from the group of left-continuous interval exchanges to the right-continuous interval exchange 
 group $\mathcal{E}$. 
  
 The automorphism $\Psi_T$ is of interest because it is not an inner automorphism. 
 One way to see this is through the homomorphism 
 $\phi\!:\! \mathcal{E} \rightarrow \mathbb{R} \wedge_\mathbb{Q} \mathbb{R},$ defined by
 \[  \phi(f_{(\pi, \lambda)} )= \sum_{i=1}^{n} \lambda_i \wedge_\mathbb{Q} \omega_i.      \]
See Arnoux {\bf\cite{Arnoux 81}} for a discussion of the properties of this map. The rotation $r_\alpha$ is defined by the data 
$\pi = (2, 1)$, $\lambda = ( 1-\alpha, \alpha)$, and it may be checked that $\phi(r_\alpha) = 1 \wedge \alpha. $
Any inner automorphism preserves $\phi$, but the action of $\Psi_T$ changes the sign of the scissors invariant; for instance, 
\[ \phi(\Psi_T(r_\alpha)) = \phi(r_{-\alpha}) = - \phi(r_\alpha). \] 

$\Psi_T$ is of further interest because it represents the only nontrivial class of outer automorphisms. Section 6 presents a proof 
of the following result. 
 
\begin{thm} \label{automorphisms}
$ \text{Aut}(\mathcal{E} ) =  \text{Inn}(\mathcal{E}) \rtimes \langle \Psi_T \rangle \cong \mathcal{E} \rtimes \mathbb{Z}/2\mathbb{Z}.$
\end{thm} 
 
Note that the inner automorphisms and the automorphism $\Psi_T$ act via conjugation by a transformation of $\mathbb{T}^1$. Thus, all 
automorphisms of $\mathcal{E}$ are geometric, in the sense that they are induced by the action of $\mathcal{E}$ on $\mathbb{T}^1$. 

%Consequently, the action of $\mathcal{E}$ in $\mathbb{T}^1$ by interval exchange transformations is  group structure of $\mathcal{E}$ 

%******Would like to discuss other corollaries and consequences of the lack of distortion elements. \\\
% 
%****** acknowledgements and open questions(mapping class groups, braid groups, free groups, polynomial, intermediate, and exponential 
%growth?? \\

%
%******* What conditions exist on the possible numbers of non-resolving fund. discontinuities????

%
%****** What do these results say about translation surfaces????? 

%%%%%%%%%%%%%%%%%%%%%%%%%%%%%%%%%%%

\section{Discontinuity Growth}    

For a map $f\in \mathcal{E}$, let $D(f)$  denote the set of points at which $f$ is discontinuous as a 
map $\mathbb{T}^1 \rightarrow \mathbb{T}^1$. Let $D_{np}(f)$ be those
discontinuities of $f$ which are not periodic:
\[ D_{np}(f) =D(f) \setminus \text{Per}(f). \]
Note that if $f$ is an infinite-order map and $D(f)$ is nonempty, then $D_{np}(f)$ is also nonempty.
%If $f$ is an infinite-order map where $\widetilde{D}(f)$ is nonempty, then
% $D_{np}(f)$ is also nonempty. To see this, first note that $\widetilde{D}(f) = D_{np}(f)$ 
% if $\text{Per}(f)$ is empty. Otherwise, the   
%set $\mathbb{T}^1 \setminus \text{Per}(f)$ of points with nonperiodic orbits is a member of 
%$\mathcal{P}$ and a proper subset of $\mathbb{T}^1$.
% Some left boundary point of this set is mapped into its interior by $f$, since there
% are only finitely many such boundary points. This point is necessarily an infinite-order
%  discontinuity point, which shows $D_{np}(f)$ is nonempty.
If $x \in D_{np}(f)$, both the forward and backward orbits of $x$ eventually consist entirely
of points at which $f$ is continuous, since $D_{np}(f)$ is a finite set 
of points with nonperiodic orbits. 
Moreover, for each $x \in D_{np}(f)$, there is some $k\geq 0$, such that $f^{-k}(x)$ is the last point of $D_{np}(f)$ encountered 
in the negative orbit of $x$. In particular, $f$ is continuous at all negative iterates $f^{-n}(x)$ for which $n > k$. 

\begin{defn} A nonperiodic discontinuity $x \in D_{np}(f)$ is a
 \emph{fundamental discontinuity} if $f$ is continuous at all negative iterates of $x$: 
 \[ \left\lbrace f^{-i}(x)   \right\rbrace_{i=1}^\infty\, \subseteq\, \mathbb{T}^1 \setminus D(f). \]
The set of fundamental discontinuities of $f$ is denoted $D_F(f)$. 
 \end{defn}
  
Thus, any point in $D_{np}(f)$ is either a fundamental discontinuity or a forward iterate of a fundamental discontinuity.
In particular, the set of fundamental discontinuities is nonempty whenever $D(f)$ is nonempty and 
$f$ has infinite order.

Let $f_{-}$ denote the left-continuous form of $f$:
\[ f_{-}(x) = \left\lbrace \begin{array}{cl}
 \displaystyle\lim_{y\rightarrow x^{-}} f(y), & \text{if $f$ is discontinuous at $x$}; \\          
 f(x), 
 & \text{otherwise.} \end{array} \right.  \]
 Similarly, $f_{+} =f$ may used to denote the original right-continuous map. 
  Observe that $(f_{-})^n = (f^n)_{-}$ and  $(f_{+})^n = (f^n)_{+}$ for all integers $n$;
 such compositions are thus denoted $f_{-}^n$ and $f_{+}^n$ without ambiguity. It follows 
 that an iterate $f^n$ is continuous at $x$ if and only if 
 \[f_{-}^n(x) =  f_{+}^n(x).  \]
The sets
\[ \left\lbrace f_{+}^n(x) \right\rbrace_{n=0}^\infty \ \ \text{and} \ \ \left\lbrace f_{-}^n(x) \right\rbrace_{n=0}^\infty \]
are called the \emph{right} and \emph{left (forward) orbits} of $x$, respectively. 

Let $x \in D_{np}(f)$ be a fundamental discontinuity. By the definition of $D_{np}(f)$, the right orbit $\{f_{+}^n(x)\}$ is 
nonperiodic. Since $x$ is fundamental, $f$ is continuous at all points in the negative orbit of $x$. Thus, the 
left and right orbits coincide for negative iterates of $f$, and it follows that the left orbit of 
$x$ is also nonperiodic. Therefore, since the set $D(f)$ is finite and the left and right forward 
orbits of $x$ are nonperiodic, both of these forward orbits eventually consist entirely of points at which 
$f$ is continuous. 

% Since $f$ is aperiodic (independently of any $\pm$ continuity convention), for any $x \in D_{np}(f)$ 
% both the right and left orbits of $x$ will eventually consist entirely of points at which $f$ is continuous.
% [Note: if $f$ has fixed points, the set $D_{np}(f)$ is different depending on whether one assumes $f$ to be right or left
%  continuous, even if $f$ is in normal form. Namely, if $f$ is right-continuous, the point
%  $0 \in S^1 = \mathbb{R}/\mathbb{Z}$ is an element of $D_{np}(f)$ and the left boundary of the interval $\text{Fix}(f)$ is not.
% This situation is reversed if $f$ is left-continuous. Note that neither 
% of these points can be fundamental discontinuities. In general, this complication may be ignored; the discontinuity growth of
% $f$ is in no way affected if $\text{Fix}(f)$ is excised from the circle and the remaining interval is reattached at its endpoints.
% In short, it suffices to have assumed that $\text{Fix}(f)$ is empty. In general, the notation $D_{np}(f)$ will refer to the set which is 
% defined when $f$ is assumed to be right continuous.  ] 

\begin{defn}
The \emph{stabilization time} of an interval exchange $f$ is the smallest positive integer $n_0,$ such that 
 $f$ is continuous at $f_{+}^n(x)$ and $f_{-}^n(x)$ for all $n\geq n_0$ and for all 
fundamental discontinuities $x$.  For a fundamental 
discontinuity $x$, if $f_{+}^{n_0}(x) = f_{-}^{n_0}(x),$ then $f_{+}^n(x) = f_{-}^n(x)$ for all $n\geq n_0$, since $f$ is continuous
 at all points in question. Such a fundamental discontinuity is said to be \emph{eventually resolving}. 
 Similarly, $f_{+}^{n_0}(x) \neq f_{-}^{n_0}(x)$ implies $f_{+}^n(x) \neq f_{-}^n(x)$ for all $n\geq n_0$; in 
this case, $x$ is said to be \emph{nonresolving}. 
\end{defn}

\begin{figure}[htb] \label{restricted rotation alpha,beta}
\begin{center}
\includegraphics[height=2in,width=5in,angle=0]{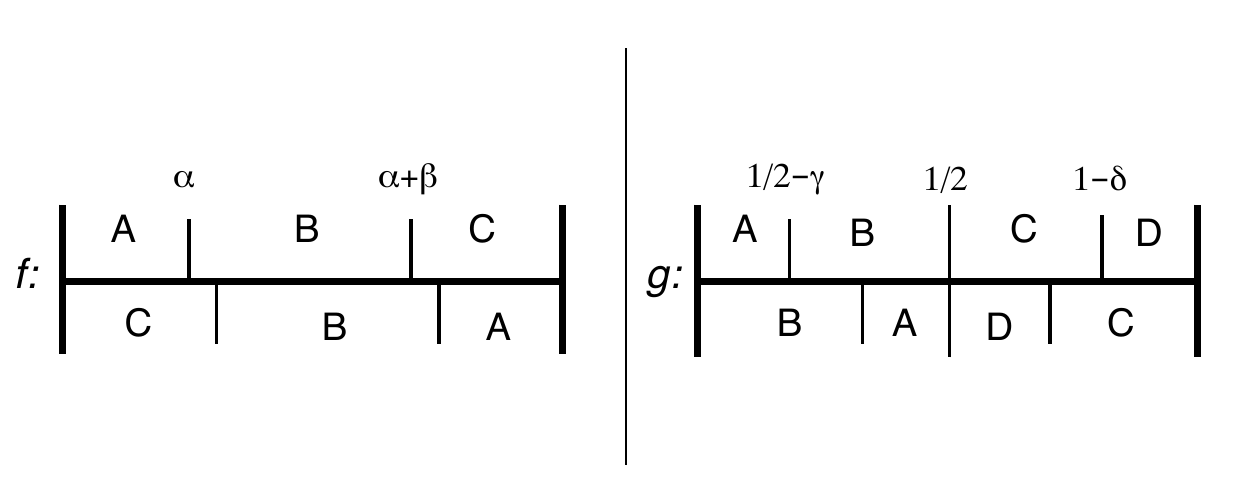}
\caption{Types of fundamental discontinuities; all parameters are irrational}
\end{center}
\end{figure}

Figure 5 gives examples of interval exchanges exhibiting the two types of fundamental discontinuity. 
 The map $f$ has 
fundamental discontinuities at $\alpha$ and $\alpha+\beta$, and it may be checked that both of these
 discontinuities are nonresolving.  The map $g$ is a product of two restricted irrational
rotations. It has fundamental discontinuities at $\frac{1}{2} - \gamma$ and $1 - \delta$, both of which are eventually resolving.

\pagebreak

Fundamental discontinuities are so-named because they completely control the asymptotics of $d(f^n)$.

 \begin{prop}     \label{discontinuity growth dichotomy: detailed version}
For any infinite-order interval exchange $f,$ exactly one of the following holds:
\begin{itemize}
\item[ (a) ] All fundamental discontinuities of $f$ are eventually resolving, in which case $d(f^n)$ is
bounded independently of $n$. 
\item[ (b) ] The map $f$ has $k \geq 1$ nonresolving fundamental discontinuities,
in which case $d(f^n)$ has linear growth on the order of $k|n|$. 
\end{itemize}

 \end{prop}

\begin{proof}  
If the map $f$ is continuous, then conclusion (\emph{a}) holds vacuously. When $D(f)$ is nonempty, the discontinuities of $f^n$ 
are contained in the set
\[ C_n = \bigcup_{i=0}^{n-1} f^{-i} \left(D(f) \right). \]
For $p\in C_n$, the left and right orbits (and hence the $f^n$-continuity status) of $p$ are determined by the left and right orbits of the
first $f$-discontinuity $x$ that $p$ meets in its forward $f$-orbit. Suppose $x = f^j(p)$, for $j\geq 0$. If $n > j$, then 
\[   f_{\pm}^n(p)) = f_{\pm}^{n-j}(f_{\pm}^{j}(p)) = f_{\pm}^{n-j}(x),\]
where $f_{\pm}^{j}(p) = x$ only because $f$ is continuous at the points 
\[ \{ p, f(p), \ldots, f^{(j-1)}(p) \}.  \]

It will first be shown that non-fundamental discontinuities of $f$ induce a uniformly bounded number of discontinuities
for an iterate $f^n$.
Let $x \in D_{np}(f)$ be a non-fundamental discontinuity. Then some
negative iterate $y = f^{-j}(x)$, $j \geq 1,$ is a fundamental discontinuity. Thus, 
for any $n \geq 1$, the discontinuity $x$ only determines the $f^n$-continuity status of points in the set 
\[ \{x, f^{-1}(x), \ldots, f^{-(j-1)}(x) \}. \] 
The status of all other preimages $f^{-k}(x)$, such that $k\geq j,$ is determined by the fundamental discontinuity $y$.
Consequently, the number of points whose $f^n$-continuity status is determined by the 
non-fundamental discontinuity $x$ is bounded  by $j$, independently of $n$. Similarly, there is a uniform 
bound to the number of points whose $f^n$-continuity status is determined by a periodic discontinuity
 of $f.$

% Consequently, to show that $d(f^n)$ is bounded when all fundamental discontinuities of $f$ are eventually
% resolving, it suffices to show that the number of $f^n$-discontinuities 
% in the set $\{ x, f^{-1}(x), \ldots, f^{-(n-1)}(x) \}$ is bounded independently of $n$, for each eventually resolving 
%fundamental discontinuity $x$.

Next, let $x$ be an eventually resolving fundamental discontinuity of $f,$ and let $n_0$ be the 
stabilization time of $f$. Then,
 \[f_{+}^n(x) = f_{-}^n(x) \]
for all $n\geq n_0$. 
Suppose $n \geq n_0$ and $k$ is such that $0 \leq k \leq (n-n_0)$. 
The right and left orbits of $f^{-k}(x)$ are 
determined by the right and left orbits of $x$: 
\[ f_{+}^n(f^{-k}x)   =    f_{+}^{n-k}(x)   =   f_{-}^{n-k}(x)   =    f_{-}^n(f^{-k}x), \]
where the middle equality holds because $n-k \geq n_0$. Thus, $f^n$ is continuous at 
$f^{-k}(x)$,  whenever  $0 \leq k \leq (n-n_0)$. 
It follows that for all $n \geq n_0$, $\{ x, f^{-1}(x), \ldots, f^{-(n-1)}(x)  \}$
contains at most $n_0$ discontinuities of $f^n$. Therefore, $d(f^n)$ is bounded if all fundamental discontinuities 
of $f$ are eventually resolving.

Alternately, suppose  that $x$ is a nonresolving fundamental discontinuity. Then
\[  f_{+}^n(x) \neq f_{-}^n(x) \]
for all $n\geq n_0$. By an argument similar to the one above, it follows that $f^n$ is discontinuous at
  $f^{-k}(x)$, for all $k$ such that $0 \leq k \leq (n-n_0)$. Thus, if $n \geq n_0$,  $f^n$ has at least $n-n_0$ 
discontinuities in the set $\{ x, f^{-1}(x), \ldots, f^{-(n-1)}(x)  \}$. Since $n_0$ is fixed relative to $n$, this 
implies that $d(f^n)$ has linear growth. 

Consequently, the presence of at least one
nonresolving fundamental discontinuity implies linear growth of $d(f^n)$, and the presence of $k$ nonresolving 
fundamental discontinuities implies $d(f^n) \sim kn$.

\end{proof}

%%%%%%%%%%%%%%%%%%%%%%%%%%%%%%%%%%%%%%%

\section{Classification of maps with bounded discontinuity growth}

The most simple example of an infinite-order interval exchange with bounded discontinuity growth is an irrational rotation 
$r_\alpha$, for 
which $d(r_\alpha^n)$ is always zero. The discontinuity growth rate of a map is invariant under conjugation, so we 
begin by stating a theorem of Li {\bf \cite{Li 99}} which gives 
necessary and sufficient conditions for an interval exchange to be conjugate to an irrational rotation.  
 For $f \in \mathcal{E}$, let $\delta(f)$ represent the number of intervals exchanged by $f$ when viewed as 
a map $[0,1) \rightarrow [0,1)$: 
\[ \delta(f) = \min \{n: f = f_{(\pi, \lambda)} \text{ for some } \pi \in \Sigma_n, \lambda \in \Lambda_n \}. \]
 
\begin{thm_nonum}[Li {\bf\cite{Li 99}}]\label{Li's theorem}
An interval exchange map $f$ is conjugate to an irrational rotation if and only if the following hold: 
\begin{itemize}
\item[ (i) ] $\delta(f^n)$ is bounded by some positive integer $N$,
\item[ (ii) ] $f^n$ is minimal for all $n \in \mathbb{N}$, and
\item[ (iii) ] There are integers $k>0$ and $M \geq 2^{N^3 + 3N^2}$ such that $\widetilde{f} = f^k$ satisfies \\
  $\delta(\widetilde{f}) = \delta( \widetilde{f}^2) =\, \cdots\, = \delta(\widetilde{f}^M) $.
\end{itemize}
\end{thm_nonum}  

The quantities $\delta(f)$ and $d(f)$ are related, but they do not differ by a uniform constant for all $f\in \mathcal{E}$. 
For a rotation $r_\alpha$, $\delta(r_\alpha) = 2$ and $d(r_\alpha)=0$, while $\delta(f) = d(f)= 3$ for any map 
$f = f_{(\pi, \lambda)}$ with permutation $\pi = (3, 2, 1)$. It may be checked that the continuity status of the
 points $0$ and $f^{-1}(0)$ account for any difference between 
$\delta(f)$ and $d(f)$; the function $\delta(f)$ always counts these points as left endpoints of a partition interval of $f$,
 but one or both of these points may fail to be 
a discontinuity of $f$ when viewed as a map $\mathbb{T}^1 \rightarrow \mathbb{T}^1$.
Consequently, some care must taken  with condition (iii) in 
restating the above theorem in terms of the discontinuity number $d$. Conceivably, one might observe $d(f^k)$ to be constant 
over a large range of $k$ while $\delta(f^k)$ is changing frequently.

This difficulty may be overcome by a good choice of the base point on $\mathbb{T}^1$.
Presenting an interval exchange as defined
on $[0,1)$ amounts to specifying a base point 0 at which to cut the circle. Choosing a new base point amounts to 
conjugation by a rotation; since the conclusion of Li's theorem is up to conjugacy, there is no loss in changing the base point. 
If $f$ is replaced with a conjugate by a rotation, it may be assumed that $f$ is continuous at all points of the orbit  $\mathcal{O}_f(0)$.
Consequently,  $d(f^n) = \delta(f^n) - 2$ for all integers $n,$ and observing $d(f^n)$ 
to be constant is now equivalent to observing that $\delta(f^n)$ is constant.

\pagebreak

\begin{thm_nonum}[Alternate Version of Li's Theorem] \label{Reformulation of Li's Theorem}
An interval exchange map $f$ is conjugate to an irrational rotation if and only if the following hold: 
\begin{itemize}
\item[ (i) ] $d(f^n)$ is bounded by some integer $N$,
\item[ (ii) ]$f^n$ is minimal for all $n \in \mathbb{N}$, and
\item[ (iii) ] after redefining the base point (conjugating by a rotation) so that $f$ is continuous on the orbit of 0, 
there are integers $k>0$ and $M \geq 2^{N^3 + 3N^2}$ such that $\widetilde{f} = f^k$ satisfies 
  $d( \widetilde{f}) = d( \widetilde{f}^2) =\, \cdots\, = d(\widetilde{f}^M) $.
\end{itemize}
\end{thm_nonum}  

Given this version of the theorem, it may now be seen to what extent the conditions \emph{(ii)} and \emph{(iii)} hold when it is only
assumed that $d(f^n)$ is bounded. To introduce some terminology, a finite union $J$ of half-open intervals is a \emph{minimal component}
of $f$ if $J$ is $f$-invariant and the $f$-orbit of any $x\in J$ is dense in $J$. It is shown in {\bf\cite{Arnoux 81}} and {\bf\cite{Novak Thesis}} 
that for any interval exchange $f$, the set of non-periodic points of $f$ decomposes into finitely many minimal components.

\begin{lem} \label{minimality_for_all_iterates}
Suppose that $f$ is minimal and $d(f^n)$ is bounded. Then for some $k \in \mathbb{N}$, all nontrivial iterates
$f^{nk}$ are minimal when restricted to each minimal component of $f^k$.
\end{lem} 

\begin{proof} Suppose that no such integer $k$ exists. 
Then $f$ is minimal, but for some \\$k_1 = m_1 > 1$, $f^{m_1}$ has multiple minimal components. Suppose that
this integer $k_1$ has been chosen to be as small as possible. Since $f$ and $f^{m_1}$ commute,
$f$ permutes the minimal components of $f^{m_1}$. This permutation induced by $f$ is
transitive since $f$ is minimal,
and it must be of order $m_1$, by the choice of $m_1$. Thus $f^{m_1}$ has exactly $m_1$ minimal components,  denoted by $J_{1,1},\, \ldots, \,J_{1,m_1}$. 

It has been assumed that no power $f^k$ is minimal for all iterates $f^{kn}$ when restricted to any of its minimal components. Thus, there
exists a smallest integer $k_2 > 1$ such that $f^{m_2}$, where $m_2 = k_1 k_2$, is not minimal
when restricted to some minimal component of $f^{m_1}$. Suppose this component is $J_{1,1}$.    
The map $f^{m_1}$ permutes the minimal components of $f^{m_2}$ which are contained in $J_{1,1}$; 
$f^{m_1}$ acts minimally on $J_{1,1}$, and so this permutation must be transitive and have order $k_2$. Additionally, the original map $f$
permutes the minimal components of $f^{m_2}$; since it also transitively permutes the minimal components of $f^{m_1},$
it follows that $f^{m_2}$ must have $k_2$ minimal components in each one of the $J_{1,j}$. Thus $f^{m_2}$ has
exactly $k_2 k_1 = m_2$ minimal components. 

By the assumption that no $k$ satisfies the conclusion of the lemma, this process may continue indefinitely. In particular, there are sequences of integers $k_i > 1$ and $m_i = \Pi_{j=1}^i k_j$, such that $f^{m_i}$ has exactly $m_i$ minimal components. 
 
 To arrive
at a contradiction with the hypothesis that $d(f^n)$ is bounded, observe that if a map $g$ has $m >1$
 minimal components $J_1,\ldots, J_m,$  then 
it must have at least $m$ discontinuities. To see this, consider a left-boundary point $x_i$ of $J_i$. 
Since some iterate of $x_i$ will eventually fall in the interior of $J_i$, it follows that 
the orbit of each $x_i$ must contain a discontinuity of $g$.
% Points just
% to the right of $x_i$ are also contained in $J_i$, but points just to the left of $x_i$ must be in a different minimal component. 
% However, after some number of iterates the orbit of $x_i$ must end up in the interior of $A_i$, in which case the points on
% both sides of this image must are in $A_i$. This is only possible if at some previous stage of its orbit $x_i$ reaches a 
% discontinuity. 
Since these orbits are distinct, the map must 
have at least $m$ discontinuities. Thus, it is impossible for $f^n$ to have an arbitrarily large number of minimal components if $d(f^n)$ is bounded. \end{proof}

\begin{rem}
It seems plausible that the above lemma should hold in general; i.e., the condition that $d(f^n)$ is perhaps not necessary. However, the argument above
strongly uses this assumption and breaks down without it. 
\end{rem}

\begin{lem} \label{constant_discontinuity_count}
Suppose $f$ has infinite order and $d(f^n)$ is bounded. Then for some $N \in \mathbb{N}$,
 $d(f^{nN})$ is constant over all $n\in \mathbb{N}$.
\end{lem}

\begin{proof} By initially replacing $f$ with an iterate, 
it may be assumed that $\text{Per}(f) = \text{Fix}(f)$. 
Let $D_F = \lbrace x_1, \ldots, x_k \rbrace $
be the fundamental discontinuities of $f$.  
 Since $d(f^n)$ is bounded, each $x_i$ is eventually resolving.
 All other non-fixed discontinuities are found in the forward orbits 
of the fundamental discontinuities.  Choose an integer $N_1 > 0$ such 
that any point of $D_{np}(f)$ may 
be reached from $D_F$ by at most $N_1$ iterates of $f$. Such an $N_1$ exists since the set $D_{np}(f)$ is finite. 

Choose $N_2$ such 
that the right and left orbits of all discontinuities in $D_{np}(f)$ are stabilized after $N_2$ iterates of $f$. 
In the situation where a non-fundamental discontinuity $x \in D_{np}(f)$ is fixed from the left (i.e., 
$f_{-}(x) = x$), it is the case that $f_{+}^n(x) \neq f_{-}^n(x)$ for all $n\geq 1$, since the right 
orbit of $x$ is nonperiodic. Otherwise, both the right and left forward orbits of any $x \in D_{np}(f)$ 
eventually consist entirely of continuity points of $f$. Thus, the notion of stabilization time 
is well-defined for all $x\in D_{np}(f)$.

Finally, choose $N > N_1 +N_2$. 
It will be shown that $d(f^{kN})$ is constant over
 all $k\in \mathbb{N}$. Since $\text{Per}(f) = \text{Fix}(f)$, the set of 
fixed discontinuities is identical for all nonzero iterates of $f$. 
Thus, it suffices to only consider the set $D_{np}(f^N)$ of
 non-fixed discontinuities of $f^N$; any such 
point must be of the form $f^{-i}(x)$, where $x \in D_{np}(f)$ and $0 \leq i < N$. The non-fixed
discontinuities of $f$ are contained in the set 
\[ \bigcup_{i=0}^{N_1} f^i(D_F). \] 
It follows that the non-fixed discontinuities of $f^N$ are contained in the set
\[ \bigcup_{i = -(N-1)}^{N_1} f^i(D_F) .\]
Let 
\[ P =  D(f^N) \cap \left( \bigcup_{i=1}^{N_1} f^i(D_F) \right),  \ \ \  Q = D(f^N) \cap \left( \bigcup_{i = -(N-1)}^0 f^i(D_F) \right).    \]
 
 Consider a point $x\in P$. Since this is a discontinuity of $f^N$, the forward $f$-orbit 
 of $x$ must encounter a discontinuity  of $f$ whose right and left orbits control
  the continuity status of $x$. Since $x$ is in $P$, this controlling discontinuity is
  non-fundamental, and it must be encountered within $N_1$ iterates of $f$.
Since $N > N_1 + N_2$, the inequality
between $f^N(x^+)$ and $f^N(x^-)$ occurs at a place where the right and left
orbits of the controlling discontinuity have already
stabilized. Thus, the right and left orbits of the controlling discontinuity are nonresolving,
and it follows that $x$ is a
discontinuity of $f^n$, for all $n \geq N$. In particular, $x$ is a discontinuity for all $f^{kN}$. Similarly, if a point 
in $ \bigcup_{i=1}^{N_1} f^i(D_F)$ is a point of continuity for $f^N$, 
it must be a point of continuity for all $f^{kN}$.
  
Next, consider a point $x\in Q$. This point is a discontinuity of $f^N$ whose $f$-orbit is controlled by a fundamental discontinuity $x_i$
of $f$. Observe that under $f^{kN}$, the image of $x$ is contained in 
\[ \bigcup_{i = (k-1)N+1}^{kN} f^i(D_F). \] 
Consequently, if $k\geq 2$, the right and left orbits of $x$ (which are controlled by the right and left orbits of the fundamental discontinuity
$x_i$) have resolved once $f^{kN}$ iterates have been applied to $x$. Thus $x,$ as well as all other points in
$\bigcup_{i = -(N-1)}^0 f^i(D_F)$, are continuity points for $f^{kN}$, $k \geq 2$. In general, the $f^{kN}$-continuity status of any point
in $ \bigcup_{i = -(kN-1)}^0 f^i(D_F) $ is controlled by the right and left orbits of a fundamental discontinuity. Since these orbits all resolve 
within $N$ iterates, it follows that 
\[ D(f^{kN}) \cap \left(  \bigcup_{i = -(kN-1)}^0 f^i(D_F) \right)  = f^{-(k-1)N}(Q). \]

The previous two paragraphs have shown that 
\[ D(f^{kN}) = P \cup f^{-(k-1)N}(Q). \]
This union is always disjoint, so the size of $D(f^{kN})$ is constant 
over all $k \in \mathbb{N}$.   Since 
\[ d(f^{kN}) = \left| D(f^{kN} ) \right| + \left| \{ \text{fixed discontinuities of } f^{kN}  \} \right|, \]
and the second term in this sum is constant over all iterates of $f$, it follows that $d(f^{kN})$ is 
constant over all $k \in \mathbb{N}$, as desired.  \end{proof}

\begin{proof}[Proof of Theorem \ref{bounded growth classification}]
 Since $f$ may be replaced with a power of itself, it
  may be assumed that $\text{Per}(f) = \text{Fix}(f)$. By applying 
 Lemma  \ref{minimality_for_all_iterates} to the restriction of $f$ on each of its minimal components,
  there is some $k$ such that any $f^{nk}$ is minimal when 
 restricted to any minimal component $J_1, \ldots, J_m$  of $f^k$.  
 Since the result is up to conjugacy in $\mathcal{E}$, it may be assumed
 that the minimal components $J_i$ are all intervals.

Consider the restriction ${f_j}$ of $f^k$ to its minimal component $J_j$. 
 It suffices to show that $f_j$ is conjugate to an irrational rotation. The function
$d( f_j^n )$ is bounded, and by construction $f_j^n$ is minimal for all $n > 0$. If necessary, conjugate $f_j$ by a 
rotation to assure that $f_j$ is continuous at all points of the orbit of 0. 
By Lemma \ref{constant_discontinuity_count}, there exists $N_j$ such that
 $d(f_j^{nN_j})$ is constant for all $n$. Consequently, the alternate version of Li's theorem 
 applies to the restricted map $f_j$, and so this map is conjugate
  to an irrational rotation.  \end{proof}

%%%%%%%%%%%%%%%%%%%%%%%%%%%%%%%%%%%%%%

\section{Proof of Theorem \ref{no distortion elements} }

We prove in this section that  $\mathcal{E}$ does not contain distortion elements.
To achieve this, it suffices to prove that an infinite-order interval exchange is not a distortion element in any finitely 
generated subgroup of $\mathcal{E}$ which contains it.

By Theorem \ref{discontinuity growth dichotomy: basic version}, 
the iterates of $f$ have linear or bounded discontinuity growth. Suppose first that $f$ 
has linear discontinuity growth.  Let $S = \{g_1,\,  \ldots,\, g_k\}$
 generate a subgroup of  $G < \mathcal{E}$ which contains $f$,  and let 
\[ M = \max_i \left \lbrace  d(g_i)   \right \rbrace. \] 
Then 
\[  d(f^n) \leq M |f^n|_{S}, \]
since $f^n$ may be expressed 
as a composition of $|f^n|_{S} $ elements from the set of generators. Consequently, linear growth of $d(f^n)$ implies linear 
growth of $|f^n|_{S},$ and thus $f$ is not a distortion element of $G.$

%In addition, for any distortion element $f$ of $\mathcal{E}$, the scissors invariant must vanish: $\phi(f) = 0$.
%Recall that $\phi$ is a group homomorphism
%\[ \phi: \mathcal{E} \rightarrow \mathbb{R} \wedge_\mathbb{Q}  \mathbb{R}. \] 
%The image of a distortion element under any group homomorphism must either be a distortion element or have
%finite order. Since, $\mathbb{R} \wedge_\mathbb{Q}  \mathbb{R}$ is a torsion-free abelian
%group, it has no distortion elements.

%Consequently, any distortion element $f$ in $\mathcal{E}$ must have bounded discontinuity growth and a vanishing scissors invariant.
%Thus, some power $f^n$ is conjugate to a product of disjointly supported irrational rotations. The map
% $f^n$ is also a distortion element, and $f^n$ also has vanishing scissors invariant. 
% These conditions seem too strong to be simultaneously satisfied; a product of irrational rotations
%will only have a vanishing scissors invariant if the rotations satisfy some numeric relation. For example, a map consisting of two disjoint
%irrational rotations which are inverses of each other will have a vanishing scissors invariant, and it seems suspicious that such a map could be a 
%distortion element while its two factors cannot be. This evidence suggests the following 

Suppose now that $f$ has infinite order and bounded discontinuity growth, and again suppose 
$f \in G = \langle g_1, \ldots, g_n \rangle < \mathcal{E}.$  
 By Theorem \ref{bounded growth classification}, after conjugation 
and replacing $f$ by an iterate it may be assumed that $f$ is a product of disjointly supported infinite-order restricted rotations. 

Let $r_{\alpha, \beta}$ denote one of these rotations, and assume first that $\alpha \notin \mathbb{Q}$. 
Let $V$ be the $\mathbb{Q}$-vector subspace of $\mathbb{R}/ \mathbb{Q}$ which is generated by 
the set of distances an
element of $G$ may translate a point of $\mathbb{T}^1$.  The space $V$ is a finite-dimensional
 $\mathbb{Q}$-vector space, since it is generated by the components of the translation vectors $\omega(g_i), 1 \leq i \leq n$.
 
Fix a basis for $V$ which includes the class $[\alpha] \in \mathbb{R}/ \mathbb{Q}$ and the class $[\beta]$ if $\beta \notin \mathbb{Q}$. Let 
\[ P_\alpha: V \rightarrow \mathbb{Q} \]
be the linear projection which returns the $[\alpha]$-coordinate of a vector with respect to this basis.
For $p \in \mathbb{T}^1,$ define the function $\phi_{\alpha,p} \!:\! G \rightarrow \mathbb{Q}$ by
\[  \phi_{\alpha,p} (g) = P_\alpha(g(p)- p). \] 

Note that the maps $\phi_{\alpha, p}$ satisfy the cocycle relation \[\phi_{\alpha, p}(fg) = \phi_{\alpha, p}(g) + \phi_{\alpha, g(p)}(f).\]
The map $f$ rotates by $\alpha \text{ mod } \beta$ on the interval $[0, \beta)$ and $P_\alpha(\beta) = 0,$  so it follows that 
\[  \phi_{\alpha,0}(f^n) = n, \ \text{for all } n \in \mathbb{Z}. \]

Now consider the generators $g_1, \ldots, g_n$. Each one of these maps induces only finitely many 
distinct translations, namely the components of $\omega(g_i)$. Consequently, there is a constant $M >0$ such that
\[   \left| \phi_{\alpha, p} (g_i) \right|   \leq  M,  \text{ for all }  p \in \mathbb{T}^1, 1 \leq i \leq n.  \] 
Thus, for any $g \in G$, 
\[ \left| \phi_{\alpha, 0} (g) \right| \leq M |g|_{S}. \] 
In particular, 
\[ n = \phi_{\alpha, 0}(f^n) \leq M |f^n|_{S}, \text{ for all } n \in \mathbb{Z}, \]
which implies linear growth for $|f^n|_{S}$. Consequently, the map $f$ is not a distortion element in $G$.

Suppose we are in the case where $f$ is a product of infinite-order rotations, but all of these rotations are 
by some $\alpha_i \in \mathbb{Q}$ (mod $\beta_i \notin \mathbb{Q})$. The argument above fails in this case because
the map $\phi_{\alpha,p}$ is not well-defined when $\alpha$ is rational. However, a similar argument 
can be made by tracking the contribution from the irrational number $\beta$. Choose a new basis for $V$ which 
contains $[\beta]$, and consider the map $\phi_{\beta,0}$. The rotation by $\alpha$ mod $\beta$ on $[0, \beta)$ 
contributes $(-1)\beta$ for every loop the iterated rotation makes around this interval.  Thus, there exists some constant
$C>0,$ (for instance, any $C > \beta / \alpha$), such that
\[  |\phi_{\beta, 0}(f^n)| \geq \frac{n}{C},   \]
It is still the case that there is a constant $M>0$ such that  
\[   \phi_{\beta, 0} (f^n)   \leq  M|f^n|_{S},   \]
which again implies linear growth for $|f^n|_S$. Thus, no infinite-order element of $\mathcal{E}$ is a distortion element.

\section{Classification of Centralizers in $\mathcal{E}$ }

\subsection{The bounded growth case}
For $f\in \mathcal{E}$, let $C(f)$ denote the centralizer of $f$ in the group $\mathcal{E}$:
\[ C(f) = C_\mathcal{E}(f) := \{g\in \mathcal{E}: fg =gf\}. \]
If $f$ is minimal, then the structure of $C(f)$ is primarily determined by the discontinuity 
growth of $f$. 
In considering the situation where $d(f^n)$ is bounded, the first case to consider is when $f = r_\alpha$
is an irrational rotation. Let $R = \lbrace r_\alpha: \alpha \in \mathbb{R}/\mathbb{Z} \rbrace$ denote the subgroup
 of rotation maps.

\begin{lem} \label{Centralizer of irrational rotation} If $\alpha$ is irrational, then $C(r_\alpha)$ 
is the rotation group $R$.  \end{lem}    

\begin{proof} (See {\bf \cite{Li 99}} for an alternate proof.)
Suppose that $g\in \mathcal{E}$ commutes with $r_\alpha,$ in which case $g = r_\alpha^{-1} g r_\alpha$.
 Since $r_\alpha$ is continuous as a map $\mathbb{T}^1 \rightarrow \mathbb{T}^1$, this conjugacy implies that 
 the discontinuity set $D(g)$ is $r_\alpha$-invariant. Consequently, if $D(g)$ is a nonempty set, it must be infinite, which 
 is impossible. Thus, $D(g)$ is empty, which implies $g\in R$, as rotations are the only continuous interval exchanges. 
\end{proof}

If $f$ is minimal and $d(f^n)$ is bounded, by Theorem \ref{bounded growth classification}
 some power $f^k$ is conjugate 
to a product of disjointly supported infinite-order restricted rotations. Suppose that 
$k$ is chosen to be as small as possible, and let $J_i$, $1\leq i \leq l,$ denote the minimal components of $f^k$.
Replace $f$ by a conjugate so that the $J_i$ are intervals, and let
$r_i$ denote the restricted rotation supported on $J_i$ induced by $f^k$. Since
$f$ is minimal and commutes with $f^k$, $f$ transitively permutes the $J_i$ and induces conjugacies
 between all of the $r_i$. Consequently, the $J_i$ are all intervals of length $1/l$, and each $r_i$ rotates 
 by the same proportion of $1/l$. Let $R_i$ denote the rotation group supported on the interval $J_i$. 
 Then $f^k$ is an element of the diagonal subgroup of 
 \[ R_1 \times \cdots \times R_l, \]     
and it follows from Lemma \ref{Centralizer of irrational rotation} that  
\[ C(f^k) = \left(R_1 \times \cdots \times R_l \right) \rtimes \Sigma_l, \]
where $\Sigma_l$ is the embedding of the symmetric group which permutes the 
$J_i$ by translation.     

Since $C(f)$ is a subgroup of the virtually abelian group
$C(f^k) \cong (\mathbb{R} / \mathbb{Z})^l \rtimes \Sigma_l $, it follows that $C(f)$ is also virtually abelian. 
In addition, $f\in C(f^k)$ implies that $f$ has the form
\[ f = r_1\cdots r_l \sigma, \]
where $r_i \in R_i$ and $\sigma$ is a permutation of the $J_i$ by translation. In particular, $f$ commutes 
with the diagonal subgroup in $R_1 \times \cdots \times R_l$, and we have proved the following.  

\begin{prop} \label{centralizer of minimal and bounded disc growth has rotation subgroup}
If $f$ is minimal and $d(f^n)$ is bounded, then $C(f)$ is virtually abelian and contains a subgroup isomorphic to $\mathbb{R}/\mathbb{Z}$. 
\end{prop}

\subsection{The linear growth case}

\indent Next, suppose that $f$ is minimal and $d(f^n)$ exhibits linear growth. 
 The discontinuity structure of $f$ and its
 iterates is significantly more complicated than the bounded case. Any map $g$ which commutes with $f$ 
 must preserve this structure, and thus one would expect the centralizer of $f$ to be significantly smaller than 
 in the bounded discontinuity situation. 
 
 \begin{prop} \label{minimality and linear discontinuity growth implies virtually cyclic centralizer}
 If $f$ is minimal and $d(f^n)$ has linear growth, then $C(f)$ is virtually cyclic.
 \end{prop} 
 
 %Should I give an example of a map for which $C(f) / <f>$ is nontrivial?  

 To prove this,  let $D = D(f)$ be the discontinuity set of $f$ 
 and let $D_{NR} = \{x_1, \ldots, x_k\}$ be the set of nonresolving fundamental discontinuities of $f$, which is nonempty by 
 Proposition \ref{discontinuity growth dichotomy: detailed version}.
  Let $n_0$ be the \emph{symmetric stabilization time} for $f$: $n_0$ is the 
 minimal positive integer such that $f$ is continuous at 
 $f^i(x)$, for all $i$ such that $|i| \geq n_0$ and all $x \in D$. The following lemma states that if a sufficiently long piece of $f$-orbit
 contains enough discontinuity points of a large power of $f$, then the $f$-orbit must contain a nonresolving fundamental
 discontinuity of $f$.

 \begin{lem} \label{Orbit with many discontinuities contains a fundamental one} 
 Suppose $f$ is minimal and has symmetric stabilization time $n_0$.
 Let $M> 3n_0$, and suppose that for some $y\in \mathbb{T}^1$ the set 
 \[ B =  \lbrace y, f^{-1}(y), \ldots, f^{-M+n_0}(y) \rbrace \]
 contains strictly more that $2n_0 +1 $ discontinuities of $f^M$. Then some $f^k(y)$, where $|k| \leq M,$ is 
 a nonresolving fundamental discontinuity of $f$.  
 \end{lem}          

 \begin{proof} Let $y_m$ denote $f^m(y)$, and
  let $j\in \mathbb{N}$ be the smallest positive integer such that $f^M$ is discontinuous at $y_{-j}$.
 Since $f^M$ is discontinuous at $y_{-j},$ this point has a controlling $f$-discontinuity at $y_{\tilde{k}} \in D(f),$
 where $-j \leq \tilde{k} \leq -j + M - 1$. Consequently, there must be a fundamental discontinuity of $f$ at 
 some $y_k$, where $-j-n_0 \leq k \leq -j + M -1$.
 
 Since $B$ contains more than $2n_0$  $f^M$-discontinuities at points $y_{-i}$ with $i >j$, there are more than $n_0$
 $f^M$-discontinuities whose status is controlled by $y_k$. In particular, at least one of the $f^M$-discontinuities 
 in $B$ is induced by the stabilized behavior of $y_k$, which implies that $y_k$ is a nonresolving fundamental 
 discontinuity of $f$. \end{proof}

\begin{proof}[Proof of Proposition \ref{minimality and linear discontinuity growth implies virtually cyclic centralizer}.]
 Suppose that $gf = fg$. Let $N\in \mathbb{N}$ be such that
 \[ N \gg n_0 \ \text{and}\ N \gg d(g).\]
 Let $x \in D_{NR}$. Since $x$ is 
 a nonresolving, the set 
 \[ A = \lbrace x, f^{-1}x, f^{-2}x, \ldots, f^{-(N-n_0)} x \rbrace \] 
 consists entirely of discontinuity points of $f^N$. 
 
 Since $f$ and $g$ commute, $g^{-1} f^N g = f^N$  is discontinuous at all points of $A$.
 Consider how this composition acts upon the set $A:$
 \begin{align*}
A = \lbrace x, f^{-1}x,& \ldots, f^{-(N-n_0)} x \rbrace  \\
&\downarrow g \\
g(A) = \lbrace gx, f^{-1}(gx),& \ldots, f^{-(N-n_0)}(gx)   \rbrace \\
  &\downarrow f^N \\
f^N g(A) = \lbrace f^N(gx), f^{N-1}&(gx), \ldots , f^{n_0}(gx)    \rbrace  \\
    &\downarrow g^{-1}\\
f^N(A) = \lbrace f^Nx, f^{N-1}&x, \ldots, f^{n_0}x  \rbrace   
\end{align*}    
 
 The cardinality of $A$ is significantly larger than $d(g)$,  which implies that $g$ acts continuously on 
 most points of $A$ in the first stage of the above composition. Similarly, $g^{-1}$ acts continuously
 on most points of $f^Ng(A)$ in the third stage. However, since $g^{-1}f^Ng$ is discontinuous
 at all points of $A$, it follows that $f^N$ is discontinuous at most of the points in 
 \[ \lbrace gx, f^{-1}(gx), \ldots, f^{-(N-n_0)}(gx)   \rbrace. \]
 By Lemma \ref{Orbit with many discontinuities contains a fundamental one}, it follows that some 
 $f$-iterate of $g(x)$ must be in $D_{NR}$. 
  
  The preceding paragraphs show that $g\in C(f)$ permutes the $f$-orbits of 
  the points in $D_{NR} = \{x_1, \ldots, x_k\}$. In particular, there is some integer $i$ and some $x_j \in D_{NR}$ such that
  \[ g(x_1) = f^i(x_j). \]
 This relation determines $g$ on the entire $f$-orbit of $x_1$:
 \[ g(f^n x_1 ) = f^n(g x_1) = f^{n+i}x_j.\]
 Since the orbit $\mathcal{O}_f(x_1)$ is dense, this relation  fully
 determines $g$. 
 
 For each $j$ such that $1\leq j \leq k$, let $h_j$ denote the unique interval exchange in $C(f)$ such that 
 \[ h_j(x_1) = x_j,\]
 if such a map exists. Then, if $g \in C(f)$ satisfies $g(x_1) = f^k(x_j)$, it follows that $g =f^k h_j $.
 In particular, $\{h_i\}$ is a set of representatives for the finite quotient group $C(f) / \langle f \rangle $, and consequently
 $C(f)$ is  virtually cyclic. \end{proof}

\subsection{Centralizers of finite order maps}

For $n\geq 2,$ the rotation $r_{1/n}$ induces a cyclic permutation of the intervals 
 \[ I_i = \left[\frac{i-1}{n}, \frac{i}{n} \right), \ \ 1\leq i \leq n.  \]
 
 Recall that the support of an interval exchange $f$ is the complement of its 
 set of fixed points. 
  For $1\leq j \leq n$, let $\mathcal{E}_{I_j}$ denote the subgroup of all interval exchanges whose 
 support is contained in $I_j$. Note that the orientation-preserving affine bijection $I_j \rightarrow[0,1)$ induces an isomorphism $\mathcal{E}_{I_j} \cong \mathcal{E}$. 
 
 Consider the following subgroups 
 in the centralizer $C(r_{1/n}).$ Let $\mathcal{E}^n_\Delta$ represent the maps in 
 $C(r_{1/n})$ which preserve the intervals $I_i$:
 \[ \mathcal{E}^n_\Delta = \lbrace g \in C(r_{1/n}): g(I_i) = I_i,\ \text{for}\  1\leq i \leq n \rbrace. \]
Note that $\mathcal{E}^n_\Delta$ is
  the diagonal subgroup of the product
 \[ \mathcal{E}_{I_1} \times \cdots \times \mathcal{E}_{I_n}, \]
as induced by the natural isomorphisms $\mathcal{E} \cong \mathcal{E}_{I_i}$. In short, a map in $\mathcal{E}^n_\Delta$ 
acts on each of the $I_j$ in the same manner, and so $\mathcal{E}^n_\Delta 
\cong \mathcal{E}_{I_j}\cong \mathcal{E}$.

 Next, let $P_n$ denote the subgroup of maps in $C(r_{1/n})$ which are invariant on $r_{1/n}$-orbits $\{ x + k/n: k = 0, 1, \ldots, n-1\}$: 
 \[ P_n = \left\lbrace g\in C(r_{1/n}): \forall x \in \mathbb{T}^1, \exists k\in \mathbb{Z}, \text{ such that }
 	g(x) = x + \frac{k}{n}\ \text{(mod 1)}   \right\rbrace . \]
 Fix $g\in P_n$, and consider a point $x=x_1 \in I_1$. Let 
 \[x_k = r_{1/n}^{k-1}(x_1) = x_1 + \frac{k-1}n,\ 2\leq k \leq n,\]
denote the other points in the $r_{1/n}$-orbit of $x$, and let $\sigma_{g,x}\in \Sigma_n$ denote the permutation that 
$g$ induces on $\{ x_i \}$: 
\[g(x_i) = x_{\sigma_{g,x}(i)}.\]
The permutation $\sigma_{g,x}$ commutes with the permutation $r: i \mapsto i+1\, \text{(mod $n$)}$, which implies 
that $\sigma_{g,x}$ must be a power of $r$. Thus, the transformation $g$ is described by a right-continuous 
(and hence piecewise constant) map 
\[ \sigma_g: I_1 \rightarrow  \langle r \rangle \cong \mathbb{Z}/ n\mathbb{Z}. \]
Conversely, any such right-continuous map $I_1 \rightarrow \mathbb{Z}/ n\mathbb{Z}$ with only finitely many 
discontinuities defines a map in $P_n$. Thus, $P_n$ is isomorphic to the abelian group of right-continuous functions 
$I_1 \rightarrow \mathbb{Z}/ n\mathbb{Z}$ having finitely many discontinuities.

\begin{prop} \label{Centralizer of rotation is Pn semidirect E}
$C(r_{1/n}) = P_n \rtimes \mathcal{E}^n_\Delta.$
\end{prop}

\begin{proof} First, suppose $g \in P_n \cap \mathcal{E}^n_\Delta.$ Then $g$ preserves the intervals 
$I_i$, which implies that $\sigma_{g,x} = id$ for all $x\in I_1$. Thus $g=id,$ and the subgroups $P_n$ and
$\mathcal{E}^n_\Delta$ have trivial intersection. 

Next, suppose $g$ is an arbitrary element of $C(r_{1/n})$. Construct $h\in P_n$ as follows. For $x=x_1\in I_1$, define 
$\{x_i\}$ as before and let
$\sigma_{h,x}$ be the permutation such that 
\[ g(x_i) \in I_{\sigma_{h,x}(i)}. \] 
Observe that $\sigma_{h,x}\in \Sigma_n$ is well-defined since $g$ maps 
an $r_{1/n}$-orbit $\{x_i\}$ to another $r_{1/n}$-orbit. Since $g$ commutes with $r_{1/n},$ the permutation 
$\sigma_{h,x}$ is a power of the permutation $r.$ Moreover, the function $x \mapsto \sigma_{h,x} \in 
\mathbb{Z}/n \mathbb{Z}[r]$ is right-continuous and has finitely many discontinuities, so it induces a map $h\in P_n$.
From its construction, $gh^{-1}$ preserves each interval $I_i$, and
 so $gh^{-1} \in \mathcal{E}^n_\Delta$. Thus $C(r_{1/n}) = P_n\cdot \mathcal{E}^n_\Delta$.  
 
 It remains to show that $P_n$ is a normal subgroup of $C(r_{1/n})$. Let $g\in P_n$ and let 
 $h\in \mathcal{E}^n_\Delta$. If $\{x_i\}$ is an $r_{1/n}$-orbit, then $h$ maps it to some other 
 $r_{1/n}$-orbit $\{y_i\}$, $g$ permutes the orbit $\{y_i\}$, and $h^{-1}$ maps $\{y_i\}$ back to 
 $\{x_i\}$. Thus $h^{-1}gh$ is invariant on $r_{1/n}$-orbits, which implies that $h^{-1}gh \in P_n$.
 Consequently, $P_n \trianglelefteq C(r_{1/n})$. 
 \end{proof}
% Then $\omega_g(x) = k/n$ and  $\omega_h(x + j/n) = \omega_h(x)$ for all $x$. Thus 
% \[ \omega_{h^{-1}gh}(x) = \omega_h(x) + \omega_g(hx) + \omega_{h^{-1}}(gh(x)) = \omega_g(hx) = k/n. \]

\begin{cor} \label{Centralizer of finite order map}
For $n\geq 2$, let $G_n = P_n \rtimes \mathcal{E}^n_\Delta$ denote the centralizer of the rotation $r_{1/n}$, and
let $G_1 = \mathcal{E}$. If $f$ is any finite-order map, then $C(f)$ is isomorphic to a finite direct product of the $G_i$.
\end{cor}

\begin{proof} Decompose $\mathbb{T}^1$ into finitely many nonempty
 \[ I_j = \text{Per}_j(f) = \{x \in \mathbb{T}^1: |\mathcal{O}(x)| = j \}. \]
 
 This decomposition is finite because an interval exchange cannot have periodic points of arbitrarily large minimal period.
      After replacing $f$ by
a conjugate, it may be assumed that the $I_j$ are intervals on which $f$ acts by a finite-order rotation. The $I_j$ are 
invariant under all $g\in C(f),$ and $C(f) \cap \mathcal{E}_{I_j}$ 
is isomorphic to $G_j$. \end{proof}

\subsection{The general situation}

Let $f$ be any interval exchange. Let $J_1, \ldots, J_k$ denote the minimal components of $f,$ let
$A = \text{Per}(f) \setminus \text{Fix}(f)$ and let $B = \text{Fix}(f)$. After replacing $f$ by a conjugate, it may be assumed 
that each of these sets is an
interval. Let $f_i$ be the restriction of $f$ to $J_i$, here defined on all of $\mathbb{T}^1$ by
\[ f_i(x) = 
\begin{cases} 
f(x) , &\text{if $x \in J_i$}\\ 
x, &\text{otherwise.}\\
\end{cases}    \]

Let $g\in C(f)$. The sets $A$ and $B$ are both $g$-invariant, but $g$ may permute the minimal components $J_i$. 
However, if $g$ maps $J_i$ onto $J_j$, then $g$ induces a conjugacy between $f_i$ and $f_j$. After replacing $f_j$ 
by a conjugate in  $\mathcal{E}_{I_j}$, it may be assumed that 
\[ f_i = \tau_{ij} f_j \tau_{ij} \]    
where $\tau_{ij}$ is the order-two map which interchanges $J_i$ and $J_j$ by translation and fixes all other points.
Replace $f$ by a further conjugate so that $f_i = \tau_{ij} f_j \tau_{ij}$ holds for all pairs $i \neq j$ such that $f_i$ and
$f_j$ are conjugate, and let $F$ be the group generated by all such $\tau_{ij}$. Note that $F$ is isomorphic to a direct product
of symmetric groups, since the relation $i \sim j \Leftrightarrow (f_i$ is conjugate to $f_j$) is an equivalence relation on 
$\{1, \ldots, k\}$.
Let $C_i = C_{\mathcal{E}_{J_i}}(f) = C(f) \cap \mathcal{E}_{J_i}$ denote the subgroup of maps in $C(f)$ with support in 
$J_i$, and let $C_A = C(f) \cap \mathcal{E}_A$. 

\begin{thm} \label{general centralizer form}
 For any $f \in \mathcal{E},$ let $A = \text{Per}(f) \setminus \text{Fix}(f)$ and let $B = \text{Fix}(f)$. Then, 
\[ C(f) \cong \left( \left( \prod_{i=1}^k C_i   \right) \rtimes F  \right) \times C_A \times \mathcal{E}_B, \]
where $F$ is a direct product of symmetric groups, where each $C_i$ is either an infinite virtually cyclic group or a subgroup of
 $(\mathbb{R} / \mathbb{Z})^n \rtimes \Sigma_n$  containing the diagonal in 
 $(\mathbb{R} / \mathbb{Z})^n$, and where $C_A$ is a direct product of finitely many factors $G_n =  P_n \rtimes \mathcal{E}^n_\Delta.$
 The factors $C_A$ and $\mathcal{E}_B$ are trivial if $\text{Per}(f) = \emptyset$, and the factors $C_i$ and $F$ are trivial if $f$ has finite order.
\end{thm}  

\begin{proof} It is clear that 
\[ C(f) \cong C_{\cup J_i}(f) \times C_A \times \mathcal{E}_B, \]
since these are disjoint and non-conjugate $f$-invariant sets which cover $\mathbb{T}^1$. The verification that 
\[ C_{\cup J_i}(f) \cong \left( \prod_{i=1}^k C_i   \right) \rtimes F \]
is similar to the proof of Proposition \ref{Centralizer of rotation is Pn semidirect E}. \end{proof}

A corollary to this result states that the existence of periodic points for an interval exchange $f$ is characterized by its centralizer.

\begin{cor} \label{periodic points iff big centralizer}
 For any $f\in \mathcal{E}$, $\text{Per}(f)$ is nonempty if and only if $C(f)$ contains a subgroup isomorphic to 
 $\mathcal{E}$.   
\end{cor}

\begin{proof} The factors $C_A$ and $\mathcal{E}_B$ both contain subgroups isomorphic to $\mathcal{E}$ if they are nontrivial, 
and at least one of these factors is nontrivial when $\text{Per}(f)$ is nonempty.  

It remains to show that if $\text{Per}(f)$ is empty, then no subgroup of $C(f)$ is
isomorphic to $\mathcal{E}$.
In this case, 
\[ C(f) \cong \prod_{i=1}^k C_i \rtimes F,  \] 
where each $C_i$ is either virtually cyclic or isomorphic to a 
subgroup of  $(\mathbb{R} / \mathbb{Z})^n \rtimes \Sigma_n$ containing the diagonal.
It may be seen that for any two infinite-order 
$g, h \in C(f)$, there are nontrivial powers $g^j$ and $h^k$ of these maps which commute. This property
does not hold for the group $\mathcal{E}.$ For instance, consider an irrational rotation $r_\alpha$ 
and any infinite-order map $f\in \mathcal{E}$ which is not a rotation; by Lemma \ref{Centralizer of irrational rotation},
nontrivial powers of $r_\alpha$ and $f$ do not commute.  Thus, it is not possible to 
embed $\mathcal{E}$ as a subgroup of $C(f)$ when $f$ has no periodic points. 
\end{proof}

\begin{cor} \label{centralizers have infinite index}
For any $f \in \mathcal{E}$ such that $f \neq id$, the index $[\mathcal{E}\!: \!C(f)]$ is uncountabe. 
\end {cor}

\begin{proof} From the structure of $C(f)$ given in the proposition, it suffices to consider the cases where $f$ has finite order, 
where $f$ is minimal with $d(f^n)$ bounded, and where $f$ is minimal with linear discontinuity growth.

If $f$ has finite order, it suffices to consider the case $f = r_{1 / n}$.  Fix an irrational $\alpha$ in $(0, 1),$ and 
note that the product of restricted rotations $r^{-1}_{\alpha \epsilon, \epsilon} r_{\alpha \epsilon', \epsilon'}$ is never an element of $C(f)$ for any 
$0 < \epsilon < \epsilon' < \frac{1}{n}$. Consequently, 
the $r_{\alpha \epsilon, \epsilon}$ provide an uncountable set of distinct coset representatives for $\mathcal{E} / C(f)$.

%If $f$ has finite order, it suffices to consider the case $f = r_{1 / n}$. For any $h \in \mathcal{E}$,  define
%\[  \alpha(h) = \min_{x \neq y:\, x = y + j/n} \{ \rho_{\mathbb{T}^1}(h(x), h(y)) \}. \]
%The function $\alpha$ measures how close together $h$ is able to map a pair of points in the same $r_{1/n}$-orbit.     
%Since transformations in $C(r_{1/n})$ preserve the collection of sets $\{x + i/n\}$, the function $\alpha$
%is constant on cosets $h \cdot C(f)$. The function $\alpha$ takes all values in $(0,1/n)$, and thus the index of 
%$C(f)$ in $\mathcal{E}$ is uncountably infinite. 

If $f$ is minimal and $d(f^n)$ is bounded, consider a conjugate $g$ of $f^k$ which is a product of 
infinite-order restricted rotations on intervals of length $1/l$.
Again, notice that for $0 < \epsilon < \epsilon' < 1/l,$ 
the product $r^{-1}_{\alpha \epsilon, \epsilon} r_{\alpha \epsilon', \epsilon'}$ is not an element of $C(g).$  Consequently, 
the $r_{\alpha \epsilon, \epsilon}$ also provide an uncountable set of coset representatives for $\mathcal{E} / C(g)$, and 
it follows that $C(f^k)$ and $C(f)$ have uncountable index in $\mathcal{E}$.

If $f$ is minimal with linear discontinuity growth,  then by proposition
 \ref{minimality and linear discontinuity growth implies virtually cyclic centralizer}, $C(f)$ is virtually cyclic. In particular,
 $C(f)$ is countable, which implies that $C(f)$ has uncountable index in $\mathcal{E}.$  \end{proof}

\section{Computation of $\text{Aut}(\mathcal{E})$}

The proof of Theorem \ref{automorphisms} is based on observing that an arbitrary 
$\Psi \in \text{Aut}(\mathcal{E})$ preserves the structure of centralizers, which implies that $\Psi$ preserves various 
dynamical properties of individual maps and subgroups in $\mathcal{E}$. 

\begin{lem} \label{irrational rotations are characterized by centralizers}
An interval exchange $f$ is conjugate to an irrational rotation $r_\alpha$ if and only if the following conditions hold:
\begin{itemize}
\item[(1) ] $C(f) \cong \mathbb{R}/\mathbb{Z} $;
\item[(2) ] if $g\in C(f)$ has infinite order, then $C(g) = C(f)$. 
\end{itemize}
\end{lem}

\begin{proof} By Lemma \ref{Centralizer of irrational rotation}, conditions (1) and (2)
hold if $f = r_\alpha$ is an irrational rotation, and these conditions are both preserved under conjugation.

Conversely, assume that $f$ satisfies (1) and (2). By
Corollary \ref{periodic points iff big centralizer}, $\text{Per}(f)$ is empty. 
Next, suppose that some $f^n$ has at least two minimal components, and denote them by $J_i$.
 Let $g$ be the map which is equal to $f$ on $J_1$ and fixes all other points. Then $g$ has infinite order and commutes
 with $f$, so $C(g) \cong \mathbb{R}/ \mathbb{Z}$ by condition (2). However, $g$ has fixed points, and so $C(g)$ 
 contains a subgroup isomorphic to $\mathcal{E}$, which is impossible by Corollary \ref{periodic points iff big centralizer}.
  Thus, $f^n$ is minimal for all $n \geq 1$. 
 
 Furthermore, $f$ has bounded discontinuity growth. If not, then $C(f)$ 
 is virtually cyclic by 
 Proposition \ref{minimality and linear discontinuity growth implies virtually cyclic centralizer},
 which is not the case for $\mathbb{R}/\mathbb{Z}$. 
 Consequently, by Theorem \ref{bounded growth classification} some power $f^k$ is conjugate to an irrational rotation.
Since $C(f) = C(f^k)$, it follows that $f$ is also conjugate to an irrational rotation. \end{proof}

Let $R < \mathcal{E}$ denote the group of circle rotations $\{r_\alpha: \alpha \in \mathbb{R}/\mathbb{Z}\}$.
For any $f\in \mathcal{E}$, let $\Phi_f$ denote conjugation by $f^{-1}$; i.e., $\Phi_f(g) = fgf^{-1}$.

\begin{cor} \label{psi maps R to a conjugate}
For any $\Psi \in \text{Aut}(\mathcal{E})$, $\Psi$ maps the rotation group $R$ to a conjugate. That is, there
exists $g\in \mathcal{E}$ such that $\Psi(R) = gRg^{-1}$.  
\end{cor}

\begin{proof} Since conditions (1) and (2) in 
Lemma \ref{irrational rotations are characterized by centralizers} are purely group theoretic, they are 
preserved by any automorphism $\Psi$. Fix an irrational rotation $r_\alpha$. By the Lemma, $\Psi(r_\alpha)$
is conjugate to an irrational rotation. In particular, there is some $g\in \mathcal{E}$ and some irrational
 $\beta \in \mathbb{R}/ \mathbb{Z}$ such that 
 \[ \Psi(r_\alpha) = \Phi_g (r_\beta). \]
Then
\[ \Psi(R) = \Psi(C(r_\alpha)) = C(\Psi(r_\alpha)) = C(\Phi_g(r_\beta))  = gRg^{-1}. \] \end{proof}

A similar result holds for maps that are conjugate to an infinite-order restricted rotation $r_{\alpha, \beta}$.

\begin{lem} \label{restricted rotation characterized by centralizers}
An interval exchange $f$ is conjugate to an infinite-order restricted rotation $r_{\alpha, \beta}$ 
if and only if the following hold: 
\begin{itemize}
\item[(1) ] $C(f) = \mathcal{E}_* \times H$, where $\mathcal{E}_* \cong \mathcal{E},$ 
     $H \cong \mathbb{R}/ \mathbb{Z}$, and $f \in  H;$ 
\item[(2) ] if $g \in H$ has infinite order, then $C(g) = C(f)$; 
\item[(3) ] for $h \in C(f),$ if the index $[C(f)\negmedspace:\negmedspace C(h) \cap C(f)]$ is finite and 
$C(h) \supsetneqq C(h) \cap C(f)$,  then $h$ is a finite-order element of $H$. 
\end{itemize}
\end{lem}

\begin{proof} Suppose that $f = r_{\alpha, \beta}$ with $\beta < 1$ and $\alpha / \beta$ irrational.
 Let $I = [\beta, 1).$ Then
\[ C(r_{\alpha, \beta}) = \mathcal{E}_I \times R_\beta, \]
where $R_\beta \cong \mathbb{R} / \mathbb{Z}$ is the group of all restricted rotations $r_{\gamma, \beta}$ on $[0, \beta)$. 
Any other infinite-order element of $R_\beta$ has the same centralizer as $r_{\alpha, \beta}$, and it follows that  
$r_{\alpha, \beta}$ satisfies conditions (1) and (2). 

To verify condition (3) for $r_{\alpha, \beta}$, take $h\in C(r_{\alpha, \beta})$ and write 
$h = h_I r_{\gamma, \beta}$, where $h_I \in \mathcal{E}_I$ and $r_{\gamma, \beta} \in R_\beta$.
Assume that $C(h)$ satisfies the hypotheses of condition (3). Note that 
\[ C(h) \cap C(r_{\alpha, \beta}) = C_{\mathcal{E}_I}(h_I) \times R_\beta, \] 
and consequently, 
\[ [C( r_{\alpha, \beta})\!:\! C(h) \cap C( r_{\alpha, \beta}) ] = [\mathcal{E}_I\!:\! C_{\mathcal{E}_I}(h_I)]. \]
Corollary 
\ref{centralizers have infinite index} states that the index   $[\mathcal{E}_I\!:\! C_{\mathcal{E}_I}(h_I)]$ is infinite if $h_I$ is not the identity. 
However, it has been assumed that this index is finite; thus
$h_I = id$ and $h = r_{\gamma, \beta}$ is a restricted rotation.
 It has also been assumed that
\[ C(h) \supsetneqq C(h)\cap C(r_{\alpha, \beta}) = C(r_{\alpha, \beta}), \]
and this is possible only if the rotation $h = r_{\gamma, \beta}$ has finite order.

Finally, observe that conditions (1)-(3) are all preserved under conjugation in $\mathcal{E}$. Consequently, they 
hold for any conjugate of $r_{\alpha, \beta}$.

Conversely, suppose that $f$ is an interval exchange satisfying (1)-(3). Since $C(f)$ contains a subgroup isomorphic
to $\mathcal{E}$,  $A = \text{Per}(f)$ is nonempty by Corollary \ref{periodic points iff big centralizer}. The map $f$ does not
have periodic points of arbitrarily large period, so $\text{Fix}(f^k) = \text{Per}(f^k) = A$ for some $k\geq 1$. 
 Since $f^k$ fixes $A$, $\mathcal{E}_A < C(f^k)$. By condition (2), $C(f^k) = C(f),$
and it follows that $f$ fixes all points in $A$. Similarly, all infinite-order $g\in H$ must fix the set 
$A$, and consequently all maps in $H$ must fix $A$. Thus, $H$
is contained in $\mathcal{E}_B$, where $B = \mathbb{T}^1 \setminus A$. 

Suppose now that $f$ has $k \geq 2$ minimal components $J_i$,
 and let $h$ be the map which equals $f$ on the component $J_1$ and fixes all other points. Then $h$ has infinite order
 and commutes with $f$. Thus, by Theorem \ref{general centralizer form},
 \begin{align*}
 C(f) =&  \left ( \left( \prod_{i=1}^k C_i \right) \rtimes F \right) \times \mathcal{E}_A,  \ \text{and}\\
 C(h) =& C_1 \times \mathcal{E}_{A \cup J_2 \cup \cdots \cup J_k},
 \end{align*}  
 where $C_i = C(f) \cap \mathcal{E}_{J_i}$ and $F$ is a finite group which permutes the $J_i$.
 In particular, $C(h) \cap C(f)$ contains $\left( \prod C_i \right) \times \mathcal{E}_A$, which has 
 finite index in $C(f)$. In addition, $C(h)$ strictly contains $C(h) \cap C(f)$ since $h$ has a larger fixed point set
 than $f$. Thus, condition (3) implies that $h$ must have finite order, which is a contradiction. 
 A similar argument may be applied to any infinite-order $g\in H$; thus, 
 all such maps have a single minimal component, namely $B$.
 
 Consider the natural isomorphism 
 \[ \mathcal{E} \rightarrow \mathcal{E}_B. \] 
Let $\widetilde{f}$ denote the preimage of $f \in \mathcal{E}_B$, and let $\widetilde{H}$ denote
the preimage of $H.$ Then all infinite-order 
$\widetilde{g} \in \widetilde{H}$ are minimal, and 
\[ C(\widetilde{g}) = C(\widetilde{f}) > \widetilde{H}, \]
which implies that all infinite-order $\widetilde{g}$ have bounded discontinuity growth. As in the proof of 
Lemma \ref{irrational rotations are characterized by centralizers}, it follows that all $\widetilde{g} 
\in \widetilde{H}$ are simultaneously conjugate to irrational rotations. 
Back in the group $\mathcal{E}_B$, this implies that $H$ is conjugate 
to a group of restricted rotations. \end{proof}

As in the earlier case, observe that the three conditions in the previous proposition are purely group-theoretic. Consequently, 
they are all preserved by any automorphism $\Psi$, which implies the following corollary.

\begin{cor} \label{psi preserves restricted rotations}
For any $\Psi \in \text{Aut}(\mathcal{E})$ and any $f$ which is conjugate to a
restricted rotation, $\Psi(f)$ is also conjugate to a restricted rotation.
\end{cor}

Let $\mathcal{P}$ denote the set algebra consisting of all finite unions of half-open intervals $[a, b) \subseteq \mathbb{T}^1$. 

\begin{prop} \label{psi preserves restricted exchange groups}
For any $\Psi \in \text{Aut}(\mathcal{E})$ and any nonempty $A \in \mathcal{P}$, 
there is a unique $B \in \mathcal{P}$ such that
$\Psi(\mathcal{E}_A) = \mathcal{E}_B$.
\end{prop}

\begin{proof} It suffices to consider $A\in \mathcal{P}$
 to be a proper, nonempty subset of $\mathbb{T}^1$. Let $g \in \mathcal{E}$ be a map that is  
  conjugate to an infinite-order restricted rotation, such that $\text{Fix}(g) = A$.
  By Corollary \ref{psi preserves restricted rotations}, 
 $\Psi(g)$ is also conjugate to a restricted rotation; let $B = \text{Fix}(\Psi(g))$. 
 
 Observe that two infinite-order restricted rotations $g$ and $h$ commute
 if and only if one of the following holds:
 \begin{itemize}
 \item[(a) ] their supports coincide and they are simultaneously conjugate to elements in some $R_\beta$; or 
 \item[(b) ] their supports are disjoint.
 \end{itemize}
 These conditions can be characterized in terms of centralizers: (a) implies that $C(g) = C(h)$, while
 (b) implies $C(g) \neq C(h)$. In particular, each condition is preserved by any automorphism of 
 $\mathcal{E}$. 
  
 Any restricted rotation with support contained in $A = \text{Fix}(g)$ commutes with $g$ and has support disjoint 
 from that of $g$. Consequently, all restricted rotations in $\mathcal{E}_A$ must map under $\Psi$ 
 to restricted rotations with support in $B = \text{Fix}(\Psi(g))$.
 The restricted rotations in $\mathcal{E}_A$ generate this subgroup; see {\bf\cite{Novak Thesis}} for a proof of this fact. 
 Therefore, the image $\Psi(\mathcal{E}_A)$ is contained in $\mathcal{E}_B$. 
 
Similarly, under $\Psi^{-1}$ all restricted rotations with support in $B$ are mapped to restricted rotations 
which commute with $g$ and have support disjoint from that of $g$. Therefore, 
$\Psi^{-1} (\mathcal{E}_B) \subseteq \mathcal{E}_A$, and it follows 
that $\Psi(\mathcal{E}_A) = \mathcal{E}_B$. \end{proof}

\subsection{Definition and properties of $\widetilde{\Psi}$}

Given an automorphism $\Psi \in \text{Aut}(\mathcal{E})$, 
Proposition \ref{psi preserves restricted exchange groups} implies that there is a well-defined transformation
\[ \widetilde{\Psi}\!:\! \mathcal{P} \to \mathcal{P}, \] 
defined by the relation 
\[ \Psi(\mathcal{E}_A) = \mathcal{E}_{\widetilde{\Psi}(A)},   \ \  A \in \mathcal{P}. \]
In particular, $\widetilde{\Psi} (\mathbb{T}^1) = \mathbb{T}^1$ and $\widetilde{\Psi}(\emptyset) = \emptyset$, 
for all $\Psi \in \text{Aut}(\mathcal{E})$. An element  $f\in \mathcal{E}$ also induces a transformation 
 $\widetilde{f}: \mathcal{P} \rightarrow \mathcal{P},$ defined by $\widetilde{f}(A) = f(A)$.

\begin{prop} \label{properties of psi-tilde}
For all $\Psi \in \text{Aut}(\mathcal{E}),$ the transformation $\widetilde{\Psi}: \mathcal{P} \rightarrow \mathcal{P}$
has the following properties: 
\begin{itemize}
\item[(1) ] $\widetilde{\Psi}$ is an automorphism of the set algebra $\mathcal{P}$.
\item[(2) ] For any $f \in \mathcal{E}$, $\widetilde{\Psi(f)} = {\widetilde{\Psi}} \widetilde{f} \, \widetilde{\Psi}^{-1} $. 
\item[(3) ] The Lebesgue measure $\mu: \mathcal{P} \rightarrow [0,1]$ is $\widetilde{\Psi}$-invariant:
$\mu(\widetilde{\Psi}(A)) = \mu(A)$. 
\end{itemize}
\end{prop}

\begin{proof} To show that $\widetilde{\Psi}$ is a set algebra automorphism, it suffices to show
 $\widetilde{\Psi}$ preserves complements, inclusion, and unions in $\mathcal{P}$.
 If $A$ and $B$ are complements in $\mathcal{P},$ then 
the centralizer in $\mathcal{E}$ of $\mathcal{E}_A$ is $\mathcal{E}_B$, and vice versa. This same relation 
holds for $\Psi(\mathcal{E}_A) = \mathcal{E}_{\widetilde{\Psi}(A)}$ and $\Psi(\mathcal{E}_B)= 
\mathcal{E}_{\widetilde{\Psi}(B)}$, which implies $\widetilde{\Psi}(A)$ and $\widetilde{\Psi}(B)$ are complements.

 For inclusion, notice that 
\[ A \subseteq B \Rightarrow \mathcal{E}_A \leq \mathcal{E}_B \Rightarrow \Psi(\mathcal{E}_A) \leq \Psi(\mathcal{E}_B)
	\Rightarrow  \]  
\[	\mathcal{E}_{\widetilde{\Psi}(A)} \leq \mathcal{E}_{\widetilde{\Psi}(B)} \Rightarrow
	 \widetilde{\Psi}(A) \subseteq \widetilde{\Psi}(B). \]
	 
To verify that $\widetilde{\Psi}$ preserves unions, let $A$ and $B$ be elements of $\mathcal{P},$ and note that  
$\widetilde{\Psi}(A) \subseteq \widetilde{\Psi}(A\cup B)$ and $\widetilde{\Psi}(B) \subseteq \widetilde{\Psi}(A\cup B)$.
Conversely, suppose that $\widetilde{\Psi}(A\cup B) \nsubseteq \widetilde{\Psi}(A) \cup \widetilde{\Psi}(B)$. To derive a 
contradiction, let 
\[ C =  	\widetilde{\Psi}(A\cup B) \setminus \left( \widetilde{\Psi}(A) \cup \widetilde{\Psi}(B) \right). \]
Then $C \in \mathcal{P}$ is nonempty, and there exists a non-identity interval exchange \\
 $f \in \mathcal{E}_C \leq \mathcal{E}_{\widetilde{\Psi}(A\cup B) }.$ 
 The map $f$ centralizes both $\mathcal{E}_{\widetilde{\Psi}(A)}$ and $\mathcal{E}_{\widetilde{\Psi}(B)}$, so 
 the map $\Psi^{-1}(f)$ centralizes $\mathcal{E}_A$ and $\mathcal{E}_B$. This implies $\Psi^{-1}(f)$ has support 
 disjoint from both $A$ and $B$. However, this is impossible since $\Psi^{-1}(f)$ is in $\mathcal{E}_{A \cup B}$. 
 Thus,   $\widetilde{\Psi}(A\cup B) \subseteq \widetilde{\Psi}(A) \cup \widetilde{\Psi}(B)$, which completes the verification 
 that $\widetilde{\Psi}$ is a set algebra automorphism.  
  
To prove property (2), recall $\Phi_f \in \text{Aut}(\mathcal{E})$ denotes conjugation by $f^{-1}$. In particular, if $\widetilde{f}$ maps the set $A$
 to the set $B$, then $\Phi_f$ induces an isomorphism from $\mathcal{E}_A$ to $\mathcal{E}_B$.
 For any $g\in \mathcal{E}$, 
 \[ \Psi \Phi_f \Psi^{-1} (g) = \Psi( f (\Psi^{-1} g) f^{-1}) = \Psi(f)\circ g \circ \Psi(f)^{-1} .   \]
Thus $\Psi \Phi_f \Psi^{-1} = \Phi_{(\Psi f)}$, and the following diagram commutes:
\[ \begin{CD}
\mathcal{E}_A  @>{\Phi_f}>>   				\mathcal{E}_B \\
@V{\Psi}VV								  	@VV{\Psi}V \\
\mathcal{E}_{\widetilde{\Psi}(A)}    @>>{\Phi_{(\Psi f)}}>  \mathcal{E}_{\widetilde{\Psi}(B)}
\end{CD}\]
Consequently, $\widetilde{\Psi(f)} = {\widetilde{\Psi}} \widetilde{f} \, \widetilde{\Psi}^{-1} $.

To prove that $\mu$ is invariant under $\widetilde{\Psi},$ it will first be shown that if $A, B \in \mathcal{P}$ are
 disjoint and $\mu(A) = \mu(B)$, then $\widetilde{\Psi}(A)$ and $\widetilde{\Psi}(B)$ are also disjoint and have equal measure. 
 For such $A$ and $B$,  let $f \in \mathcal{E}$ be any 
interval exchange such that $\widetilde{f}(A) = B$. Then $f$ induces 
a conjugacy between the subgroups $\mathcal{E}_A$ and $\mathcal{E}_B,$ and $\Psi(f)$ induces a 
conjugacy between $\mathcal{E}_{\widetilde{\Psi}(A)}$ and $\mathcal{E}_{\widetilde{\Psi}(B)}.$   By (2),
 $\widetilde{\Psi(f)}$ maps $\widetilde{\Psi}(A)$ to $\widetilde{\Psi}(B)$, and as a result,
  $\mu(\widetilde{\Psi}(A)) = 	\mu(\widetilde{\Psi}(B))$, which proves the initial claim. 
  
  To prove that $\mu(\widetilde{\Psi}(A)) = \mu(A)$ for any $A\in \mathcal{P}$, assume first that 
  $\mu(A)$ is rational. Since any $\widetilde{\Psi}$ preserves finite disjoint unions, it may be 
  further assumed that $\mu(A) = 1/n$. Lebesgue measure is invariant under any conjugacy $\Phi_f$, so
   it finally suffices to consider the case $A = [0, 1/n).$ Each of the intervals 
  \[ A_i = \left[ \frac{i-1}{n}, \frac{i}{n} \right), \ \ 2\leq i \leq n, \]
  has the same measure as $A = A_1$ and is disjoint from it. Thus 
  \[ \mu(\widetilde{\Psi}(A_i)) = \mu(\widetilde{\Psi}(A)), \ \ 2\leq i \leq n. \]
  Since the sets $\widetilde{\Psi}(A_i) $ are also pairwise disjoint and cover $\mathbb{T}^1$, it follows 
  that $\mu( \widetilde{\Psi}(A_i)) = 1/n$. Consequently, $\widetilde{\Psi}$ preserves the measure of sets 
  with rational measure. In general, the set $A$ may be approximated by an increasing family 
  of sets in $\mathcal{P}$ having rational measure. \end{proof}

\subsection{Proof of Theorem \ref{automorphisms}}

Let $\Psi$ be an arbitrary automorphism of $\mathcal{E}$.
It will be shown that the identity automorphism may be reached by successively replacing $\Psi$
with a composition of $\Psi$ and an automorphism in   $\langle \text{Inn}(\mathcal{E}), \Psi_T \rangle$. 

To begin, by Corollary \ref{psi maps R to a conjugate}, $\Psi$ maps the rotation group $R$ to 
a conjugate $\Phi_g(R)$. Replacing $\Psi$ by $\Phi_g^{-1} \circ \Psi $, 
it may be assumed that $R$ is invariant under $\Psi$.
Let $\Psi_R : \mathbb{R}/ \mathbb{Z} \rightarrow \mathbb{R}/ \mathbb{Z} $ denote the 
restriction $\Psi |_R,$ where $r_\alpha \mapsto \alpha$ is the natural identification of $R$ and 
$\mathbb{R}/ \mathbb{Z}$. 

\begin{lem} \label{psi_R is continuous}
$\Psi_R$ is continuous (w.r.t. the standard topology on $\mathbb{R}/ \mathbb{Z}$).
\end{lem}

\begin{proof} Since $\Psi_R$ is a group homomorphism, it suffices to show that $\Psi_R$ is continuous at $0 \in \mathbb{R}/ \mathbb{Z}$.
Suppose that $\alpha_n \rightarrow 0$ in $\mathbb{R}/\mathbb{Z}$. Then for any nonempty $A \in \mathcal{P}$, there exists a 
constant $M_A > 0$ such that 
\[ A \cap r_{\alpha_n}(A) \neq \emptyset,\  \text{ if }\, n \geq M_A. \]

Conversely, this condition characterizes sequences in $\mathbb{R}/ \mathbb{Z}$ which converge to 0.
In particular, given some sequence $\alpha_n$, suppose that there exists a constant $M_A$ as above
for every nonempty $A \in \mathcal{P}$. For any $\epsilon > 0$, let $A_\epsilon = [0, \epsilon)$. Then 
\[ A_\epsilon \cap  r_{\alpha_n}(A_\epsilon) \neq \emptyset,\   \text{ if }\, n \geq M_{A_\epsilon}, \]
which implies that $|\alpha_n| < \epsilon$ for all $n \geq M_{A_\epsilon}$. Thus, $\alpha_n \rightarrow 0$.

Assuming again that $\alpha_n \rightarrow 0$, 
define $\beta_n = \Psi_R (\alpha_n)$, so $r_{\beta_n} = \Psi(r_{\alpha_n})$. Let $B \in \mathcal{P}$ be nonempty, 
and let $A = \widetilde{\Psi}^{-1}(B)$. Then by Proposition \ref{properties of psi-tilde}, part (2), 
\[ r_{\beta_n}(B) = \widetilde{\Psi}(r_\alpha(\widetilde{\Psi}^{-1}(B))) = \widetilde{\Psi}(r_{\alpha_n}(A)). \] 
Consequently, $A \cap r_{\alpha_n}(A) \neq \emptyset$ if and only if $B \cap r_{\beta_n}(B) \neq \emptyset$.
Therefore, if $\alpha_n \rightarrow 0$, then there exists $M_B$ (namely, the $M_A$ associated with
$\alpha_n$), such that 
\[ B \cap r_{\beta_n}(B) \neq \emptyset,\ \forall \, n \geq M_B. \]
From the above characterization of sequences converging to zero, 
it follows that $\Psi_R$ is continuous at zero. \end{proof}

The only continuous automorphisms of $\mathbb{R}/ \mathbb{Z}$ are the identity and $x \mapsto -x.$
Note that the restriction of the orientation-reversing automorphism $\Psi_T$ induces the second of these automorphisms. 
Subsequently, after replacing $\Psi$ by $\Psi \circ \Psi_T$  if $\Psi_R$ is not the identity, it
may be assumed that $\Psi_R = id$. 

\begin{lem} \label{if psi fixes rotations, then psi-tilde preserves intervals}
If $\Psi \in \text{Aut}(\mathcal{E})$ fixes the rotation group $R$, then $\widetilde{\Psi}$ maps
any interval in $\mathcal{P}$ to another interval. 
\end{lem}

\begin{proof} Since any rotation will preserve intervals in $\mathcal{P}$, it suffices to consider 
$I_a = [0, a)$. Then there exists some $\epsilon > 0$, such that for any $\alpha \in (-\epsilon, \epsilon),$
\[ \mu(I_a \cap r_\alpha(I_a)) = a - |\alpha|. \]
By Proposition \ref{properties of psi-tilde}, part (2), and the hypothesis that $\Psi(r_\alpha) = r_\alpha,$ 
 \[ \widetilde{\Psi} \circ \widetilde{r}_\alpha  = \widetilde{\Psi(r_\alpha)} \circ \widetilde{\Psi} = 
\widetilde{r}_\alpha \circ \widetilde{\Psi}.\]
Therefore, it is also the case that
\[ \mu( \widetilde{\Psi}(I_a) \cap r_\alpha( \widetilde{\Psi}(I_a))) = a - |\alpha|,\]
for $\alpha \in (-\epsilon, \epsilon).$ 

Suppose that $\widetilde{\Psi}(I_a)$ has $k\geq 1$ components:
 \[ \widetilde{\Psi}(I_a) = A_1 \cup \cdots \cup A_k, \]
 where the $A_i$ are pairwise disjoint intervals. Since the $A_i$ are disjoint, there is some $\delta > 0$ such
 that 
 \[ r_\beta(A_i) \cap A_j = \emptyset, \ \ \text{for all } |\beta | < \delta \text{ and } i \neq j.\]
 Consequently, if $|\beta| < \delta$, then 
 \[ \mu( \widetilde{\Psi}(I_a) \cap r_\beta(\widetilde{\Psi}(I_a))) = a - k|\beta|. \]
 It follows that $k=1$, which implies that $\widetilde{\Psi}(I_a)$ must be an interval. \end{proof}
 
 Continue with the assumption that $\Psi_R$ is the identity. By the previous lemma, $\widetilde{\Psi}$
 maps the interval $I_a = [0,a)$ to some translate of $I_a$. After composing $\Psi$ with a suitable $\Phi_{r_\beta}$,
 it may be assumed that $\Psi_R$ is the identity and $\widetilde{\Psi}(I_a) = I_a$. Since 
 \[   \widetilde{\Psi} \circ \widetilde{r}_\beta  = \widetilde{r}_\beta \circ \widetilde{\Psi}, \]
 for all $\beta \in \mathbb{R} / \mathbb{Z}$, it follows that $\widetilde{\Psi}$ fixes any translate 
 $\widetilde{r}_\beta(I_a) = [\beta, a + \beta)$. Thus, for any $\beta,\, 0< \beta < a$, 
 $\widetilde{\Psi}$ fixes the intersection 
 \[ I_a \cap r_\beta(I_a) = [\beta, a). \]
 Thus $\widetilde{\Psi}$ fixes all translates of arbitrarily small intervals, which implies that $\widetilde{\Psi}$ is
 the identity on $\mathcal{P}$. Consequently, for any $f \in \mathcal{E}$, 
 $\Psi(f)$ acts on the sets in $\mathcal{P}$ identically to the 
 way $f$ does, which implies $\Psi$ is the identity. It has been shown that any $\Psi \in \text{Aut}(\mathcal{E})$
 is in the group $\langle \text{Inn}(\mathcal{E}), \Psi_T \rangle$, and the proof of 
 Theorem \ref{automorphisms} is complete.

%:biblio

\end{document}